\documentclass[11pt]{article}

\usepackage{color}
\usepackage{amsmath, amsfonts,amssymb,theorem,euscript,array,enumerate,amsfonts,mathrsfs}
\usepackage{hhline}
\usepackage[dvips]{graphicx}
\usepackage[dvips]{graphics}
\usepackage[dvips]{epsfig}

\theoremheaderfont{\bfseries} \theoremstyle{plain}
\theorembodyfont{\upshape}
\newcommand{\esp}{\mathbb{E}} 
\newcommand{\proba}{\mathbb{P}}

\newcommand{\R}{\mathbb{R}}

\newcommand{\F}{\mathcal{F}}

\def \ep{\hbox{ }\hfill$\Box$}

\def \Frac{\displaystyle\frac}

\def\Dpi#1{\Frac{\partial #1}{\partial \pi}}

\def \Rr{\mathbb{R}}

\def \Ee{\mathbb{E}}
\def \Ff{\mathbb{F}}
\def \Pp{\mathbb{P}}
\def \Qq{\mathbb{Q}}
\def \Dd{\mathbb{D}}
\def \Gg{\mathbb{G}}

\def \Ac{{\cal A}}
\def \Bc{{\cal B}}

\def \Dc{{\cal D}}
\def \Ec{{\cal E}}
\def \Fc{{\cal F}}
\def \Gc{{\cal G}}

\def \Mc{{\cal M}}

\def \eps{\varepsilon}

\theoremstyle{plain}

\theorembodyfont{\upshape\it}
\newtheorem{Thm}{\bf Theorem}[section]

\newtheorem{Lem}{\bf Lemma}[section]

\newtheorem{Rem}{\bf Remark}[section]

\newtheorem{Exe}{\bf Example}[section]

\def\esssup_#1{\underset{#1}{\mathrm{ess\,sup\, }}}
\def\essinf_#1{\underset{#1}{\mathrm{ess\,inf\, }}}
\def\argmax_#1{\underset{#1}{\mathrm{arg\,max\, }}}

\def\reff#1{{\rm(\ref{#1})}}

\def\beqs{\begin{eqnarray*}}
\def\enqs{\end{eqnarray*}}
\def\beq{\begin{eqnarray}}
\def\enq{\end{eqnarray}}

\begin{document}

\title{Optimal investment with counterparty risk:  \\ a default-density modeling  approach}

\author{Ying JIAO
             \\\small  Laboratoire de Probabilit\'es et
             \\\small  Mod\`eles Al\'eatoires
             \\\small  CNRS, UMR 7599
             \\\small  Universit\'e Paris 7
             \\\small  e-mail: jiao@math.jussieu.fr
             \and
             Huy\^en PHAM
             \\\small  Laboratoire de Probabilit\'es et
             \\\small  Mod\`eles Al\'eatoires
             \\\small  CNRS, UMR 7599
             \\\small  Universit\'e Paris 7
             \\\small  e-mail: pham@math.jussieu.fr
             \\\small  and Institut Universitaire de France
             }

\date{\today}

\maketitle

\begin{abstract}
We consider a financial market with a stock  exposed  to a  counterparty risk indu\-cing a drop in the  price, and which can still be traded after 
this  default  time.  We use  a default-density modeling approach, and address in this incomplete market context  the expected utility maximization from terminal wealth.  We show how this problem can be suitably decomposed in  two  optimization problems in complete market framework: an after-default utility maximization and a global before-default optimization problem involving the former one.  These two optimization problems are solved explicitly,  respectively by duality and dynamic programming approaches, and provide a fine 
understanding of the optimal strategy.  We give some numerical results   illustrating the impact of  counterparty risk  and the loss given default  on optimal trading strategies, in particular with respect to the Merton portfolio selection problem. 
\end{abstract}

\vspace{13mm}





\newpage

\section{Introduction }
In a financial market, the default of a firm  has usually important influences on the other ones. This has been shown clearly by several recent default events. The impact  of a counterparty default may arise  in various contexts. In terms of  credit spreads, one observes in general  a positive ``jump'' of the  default  intensity, called the contagious jump and investigated firstly by Jarrow and Yu \cite{Jarrow-Yu}. In terms of asset (or stock) values for a firm, the default of a counterparty will in general induce  a drop of its value process. In this paper, we analyze the impact of this risk on the optimal investment problem. More precisely, we consider an agent, who invests in  a  risky asset exposed to a counterparty risk,  and we are interested in  the optimal trading  strategy and the value function when taking into account the possibility of default of a counterparty, together with the instantaneous loss of the asset at the default time.

The global market information containing default is modeled by the progressive enlargement of a background filtration, denoted by $\mathbb F$, representing the default-free information. The default time $\tau$ is in general a totally inaccessible stopping time with respect to the enlarged filtration, but is not an $\mathbb F$-stopping time. We shall work with a density hypothesis of the conditional law of default given $\mathbb F$. This hypothesis has been introduced by Jacod \cite{Jacod87} in the initial enlargement of filtrations, and has been adopted recently by  El Karoui et al. \cite{Elkaroui-Jeanblanc-Jiao} in the progressive enlargement setting for the credit risk analysis. The density approach is particularly suitable to study what goes on after the default, i.e., on $\{\tau\leq t\}$. For the before-default analysis on $\{\tau>t\}$, there exists an explicit relationship between the density approach and the widely used intensity approach.

The market model considered  here is incomplete due to the jump induced by the default time. 
The general optimal investment problem in an  incomplete market has been  studied by Kramkov and Schachermayer \cite{krasch99} by duality methods.  Recently, Lim and Quenez \cite{Lim-Quenez}  addressed, by using dynamic programming,  the utility maximization  in a market with default.     
The key idea of our paper is to derive,  by relying on  the conditional density approach of default, a natural  separation of  the initial optimization problem into an after-default one and a global before-default one. Both problems are reduced to a complete market setting, and the solution of the latter one depends on the solution of the former one. These two optimization problems are solved  by  duality and  dynamic programming approaches, and the main advantage is to give  a better insight, and more explicit results than the incomplete market framework. The interesting feature  of our decomposition is to  provide a nice interpretation of optimal strategy switching at the default time $\tau$. 
Moreover, the explicit solution (for the CRRA utility function) makes clear the roles played by the default time $\tau$ and the loss given default in the investment strategy, as shown by some numerical examples.

The outline  of this article  is organized as follows. In Section 2, we present the model and the investment problem,  and introduce  the default density hypothesis. We then explain in Section 3 how to decompose the optimal investment problem into the before-default and after-default ones.  We solve these  two optimization problems in Section 4, by using the duality approach for the after-default one and the dynamic programming approach for the global before-default one. We examine more in detail the popular case of  CRRA utility function and finally, numerical results  illustrate the impact of counterparty risk on optimal trading strategies, in particular with respect to the classical Merton portfolio selection problem.

\section{The conditional density model for counterparty risk}

\setcounter{equation}{0}
\setcounter{Thm}{0}
\setcounter{Pro}{0}
\setcounter{Cor}{0}
\setcounter{Lem}{0}
\setcounter{Def}{0}
\setcounter{Rem}{0}
\setcounter{Exe}{0}

We consider a financial market model with a riskless bond assumed for simplicity equal to one, and a stock  subject to a counterparty risk:
the dynamics of the risky asset is affected by another firm,  the counterparty, which may default,  inducing consequently  a drop in  the asset  price.
However, this  stock still exists and can be traded after the default of the counterparty.

Let us fix a probability space $(\Omega,\Gc,\Pp)$ equipped with a brownian motion $W$ $=$ $(W_t)_{t\in [0,T]}$ over a finite horizon $T$ $<$ $\infty$,
and denote by $\Ff$ $=$ $(\Fc_t)_{t\in [0,T]}$ the natural filtration of $W$.  We are given a nonnegative and finite random variable $\tau$, representing the default time, on $(\Omega,\Gc,\Pp)$.  Before  the default time $\tau$, the filtration $\Ff$ represents the information accessible to the investors. When the default occurs, the investors add this new information $\tau$ to the reference filtration $\Ff$. We then introduce
$D_t$ $=$ $1_{\tau\leq t}$, $0\leq t\leq  T$,  $\Dd$ $=$ $(\Dc_t)_{t\in [0,T]}$ the filtration generated by this jump process, and $\Gg$ $=$
$(\Gc_t)_{t\in [0,T]}$ the enlarged progressive filtration $\Ff\vee\Dd$, representing the structure of information available for the investors over $[0,T]$.

The stock price process is governed by the dynamics:
\beq \label{edsS}
dS_t &=& S_{t-} \big(\mu_t dt + \sigma_tdW_t -   \gamma_t dD_t\big), \;\;\;\;\; 0 \leq t \leq T,
\enq
where $\mu$, $\sigma$ and $\gamma$ are $\Gg$-predictable processes.  At this stage, without any further condition on the default time $\tau$, we do not know yet that $W$ is a $\Gg$-semimartingale (see Remark \ref{remdenshyp}), and the meaning of  the sde \reff{edsS} is the following. Recall 
(cf. \cite{Mansuy-Yor}) that any $\Gg$-predictable process $\varphi$ can be written in the form: $\varphi_t$ $=$ $\varphi_t^{\Ff} 1_{t\leq\tau}$ $+$
$\varphi_t^d(\tau) 1_{t>\tau}$,  $0\leq t\leq T$, where $\varphi^{\Ff}$ is  $\Ff$-adapted, and  
$\varphi_t^d(\theta)$ is measurable w.r.t.  $\Fc_t\otimes \mathcal B(\R_+)$, for all $t$ $\in$ $[0,T]$.  
The dynamics \reff{edsS} is then written as:
\beq
dS_t &=& S_t\big( \mu_t^{\Ff} dt + \sigma_t^{\Ff} dW_t\big), \;\;\;\;\; 0 \leq t < \tau,  \label{Sbefore}\\
S_{\tau} &=& S_{\tau^-}(1- \gamma_\tau^{\Ff}),  \label{Sat} \\
dS_t &=& S_t\big(\mu_t^{d}(\tau)  dt + \sigma_t^{d}(\tau) dW_t\big), \;\;\;\;\;  \tau <  t \leq T,  \label{Safter}
\enq
where $\mu^{\Ff}$, $\sigma^{\Ff}$, $\gamma^{\Ff}$ are  $\Ff$-adapted processes,  and $(\omega,\theta)$ $\rightarrow$ $\mu_t^d(\theta)$,
$\sigma_t^d(\theta)$ are  $\Fc_t\otimes\mathcal B([0,t))$-measurables functions for all $t$ $\in$ $[0,T]$.
The nonnegative process $\gamma$ represents the (proportional) loss on the stock price induced by the default of the counterparty, and we may  assume that $\gamma$ is a stopped process at $\tau$, i.e.  $\gamma_t$ $=$ $\gamma_{t\wedge\tau}$.
 By misuse of notation, we shall thus identify $\gamma$ in \reff{edsS} with the $\Ff$-adapted process $\gamma^{\Ff}$ in \reff{Sat}.  When the counterparty defaults, the drift and diffusion coefficients $(\mu,\sigma)$ of the stock price switch from $(\mu^{\Ff}$,$\sigma^{\Ff})$  to  $(\mu^d(\tau),\sigma^d(\tau))$,
 and the after-default coefficients may depend on the default time $\tau$. We  assume that $\sigma_t$ $>$ $0$, $0\leq t\leq T$,   the following  integrability conditions are satisfied:
\beq \label{integmu}
\int_0^T \Big|\frac{\mu_t}{\sigma_t}\Big|^2 dt + \int_0^T |\sigma_t|^2  dt  & < & \infty, \;\;\;\; a.s.
\enq
and
\beq \label{gam1}
0 \leq \gamma_t & < & 1, \;\;\; 0 \leq t\leq T, \;\;\; a.s.
\enq
which ensure that the dynamics of the asset price process  is well-defined, and the stock price remains (strictly) positive over $[0,T]$ (once the initial stock price $S_0$ $>$ $0$), and locally bounded.

Consider now an investor who can trade continuously in this financial market by holding a positive wealth at any time. This is mathematically quantified by a $\Gg$-predictable process $\pi$ $=$ $(\pi_t)_{t\in [0,T]}$, called trading strategy and
representing  the proportion of wealth  invested in the stock,  and the associated wealth process $X$ with dynamics:
\beq \label{dynX}
 dX_t &=& \pi_t  X_{t^-} \frac{dS_t}{S_{t-}}, \;\;\;\;\; 0 \leq t\leq T.
\enq
By writing the $\Gg$-predictable process $\pi$ in the form: $\pi_t$ $=$ $\pi_t^{\Ff}1_{t\leq\tau}$ $+$ $\pi_t^d(\tau)1_{t>\tau}$, $0\leq t\leq T$, where
$\pi^{\Ff}$ is  $\Ff$-adapted and  $\pi_t^d(\theta)$ is $\Fc_t\otimes\mathcal B([0,t))$-measurable, and in view of \reff{Sbefore}-\reff{Sat}-\reff{Safter},
the wealth process evolves as
\beq
dX_t &=& X_t \pi_t^{\Ff} \big( \mu_t^{\Ff} dt + \sigma_t^{\Ff} dW_t\big), \;\;\;\;\; 0 \leq t <  \tau,   \label{Xbefore} \\
X_\tau &=& X_{\tau^-}(1 -  \pi_\tau^{\Ff} \gamma_\tau) \label{Xat} \\
dX_t &=& X_t \pi_t^d(\tau) \big(  \mu_t^d(\tau) dt + \sigma_t^d(\tau) dW_t \big), \;\;\;\;\; \tau < t \leq T. \label{Xafter}
\enq
We say that a trading strategy  $\pi$ is admissible, and we denote $\pi$ $\in$ $\Ac$,  if
\beqs
\int_0^T |\pi_t\sigma_t|^2 dt  \; < \;  \infty,  & \mbox{ and } & \pi_\tau \gamma_\tau <  1 \;\;\; a.s.
\enqs
This ensures that the dynamics of the wealth process is well-defined with a positive wealth at any time (once starting from a positive initial capital
$X_0$ $>$ $0$).

\vspace{2mm}

In the sequel, we shall  make the standing assumption, called {\it density} hypothesis, on the default time of the counterparty.  For any $t$ $\in$ $[0,T]$,
the conditional  distribution of $\tau$ given $\Fc_t$ admits a density with respect to the Lebesgue measure,  i.e.  there exists a
family of $\F_t\otimes\mathcal B(\Rr_+)$-measurable positive functions  $(\omega,\theta)$ $\rightarrow$ $\alpha_t(\theta)$
such that:
\beqs
\hspace{-3cm} {\bf (DH)}  \hspace{4cm} \proba[\tau\in d\theta|\F_t] &=&  \alpha_t(\theta)d\theta, \;\;\;\;\;  t \in [0,T].
\enqs
We note that for any $\theta\geq 0$, the process $\{\alpha_t(\theta),0\leq t\leq T\}$ is a $(\Pp,\Ff)$-martingale.

\begin{Rem} \label{remdenshyp}
{\rm  Such a hypothesis is usual  in the theory of initial enlargement of filtration, and was introduced by Jacod \cite{Jacod87}.  The 
{\bf (DH)} Hypothesis was  recently adopted  by  El Karoui et al. \cite{Elkaroui-Jeanblanc-Jiao} in the progressive enlargement of filtration  for credit risk modeling. 
Notice that in the particular case where the family of densities satisfies 
$\alpha_T(t)$ $=$ $\alpha_t(t)$ for all $0\leq  t\leq T$, we have $\proba[\tau>t|\F_t]$ $=$ $\proba[\tau>t|\F_{T}]$. This corresponds to the so-called 
immersion hypothesis (or the {\bf H}-hypothesis), which is a familiar condition in credit risk analysis, and means equivalently that  
any square-integrable $\mathbb F$-martingale is a square-integrable $\mathbb G$-martingale.  
The {\bf H}-hypothesis appears natural for the analysis on before-default  events  when $t < \tau$, but is actually  restrictive when it concerns after-default events on $\{t\geq\tau\}$, see  \cite{Elkaroui-Jeanblanc-Jiao} for a more detailed discussion. By considering here the whole family 
$\{\alpha_t(\theta), t \in [0,T], \theta \in \R_+\}$, we obtain additional information  for the analysis of after-default events, which is crucial for our purpose.

Let us also mention that  the classical intensity of default can be expressed in an explicit way by means of  the density. Indeed, 
the $(\Pp,\Gg)$-predictable compensator of  $D_t$ $=$ $1_{\tau\leq t}$ is given by $\int_0^{t\wedge \tau}\alpha_\theta(\theta)/G_\theta d\theta$, 
where $G_t$ $=$ $\proba[\tau>t|\F_t]$ is the conditional survival probability.  In other words, the process $M_t$ $=$ 
$D_t - \int_0^{t\wedge\tau}\alpha_\theta(\theta)/G_\theta d\theta$ is a $(\Pp,\Gg)$-martingale. 
Thus, by observing from the martingale property of 
$\{\alpha_t(\theta),0\leq t\leq T\}$ that  $G_t$ $=$  $\int_t^{\infty}\alpha_t(\theta)d\theta$ $=$ $\int_t^{\infty}\esp[\alpha_\theta(\theta)|\F_t]d\theta$, we 
recover  completely the intensity process  $\lambda_t^{\mathbb G}$ $=$ $1_{t\leq\tau}\,\alpha_t(t)/G_t$ from  the knowledge of  the process 
$\{\alpha_t(t),t\geq 0\}$. 
However, given the intensity $\lambda^{\mathbb G}$, we can only obtain some part of the density family, namely  $\alpha_t(\theta)$  for $\theta\geq t$. 

Under {\bf (DH)} Hypothesis, a $(\proba,\Ff)$-brownian motion $W$ is  a $\Gg$-semimartingale and  admits an explicit decomposition in terms of the density $\alpha$ given by (see \cite{Mansuy-Yor}, \cite{Jeanblanc-LeCam}, \cite{Elkaroui-Jeanblanc-Jiao}):
\beqs
W_t &=& \hat W_t^{\Gg} + \int_0^{t\wedge\tau}\frac{d\left<W_s,G_s\right>}{G_s}
+\int_\tau^t\frac{d\left<W_s,\alpha_s(\tau)\right>}{\alpha_s(\tau)} \; =: \;  \hat W_t^{\Gg} + A_t, \;\;\; 0 \leq t\leq T, 
\enqs
where $\hat W^{\Gg}$ is a $(\Pp,\Gg)$-brownian motion, and  $A$ is a finite variation $\Gg$-adapted process.  
Moreover, by the It\^o martingale representation theorem for brownian filtration $\Ff$,  $A_t$ is written in the form $A_t$ $=$ 
$\int_0^t a_sds$ for some $\Gg$-adapted process $a$ $=$ $(a_t)_{t\in [0,T]}$.   Let us then define the $\Gg$-adapted process
\beqs
\beta_t  &=& \frac{\mu_t+ \sigma_t a_t -  \gamma_t \lambda_t^{\Gg}}{\sigma_t}, \;\;\; 0 \leq t\leq T,
\enqs
and consider the Dol\'eans-Dade exponential local martingale: $Z_t^{\Gg}$ $=$ $\Ec(-\int \beta d\hat W^{\Gg})_t$, $0 \leq t\leq T$. By assuming that 
$Z^{\Gg}$ is a $(\Pp,\Gg)$-martingale (which is satisfied  e.g. under the Novikov criterion: $\Ee[\exp(\int_0^T \frac{1}{2}|\beta_t|^2 dt)]$ $<$ $\infty$),  this 
defines a probability measure $\Qq$ equivalent to $\Pp$ on $(\Omega,\Gc_T)$ with  Radon-Nikodym density: 
\beqs
\frac{d\Qq}{d\Pp} &=&  Z_T^{\Gg} \; = \;  \exp\Big(-\int_0^T  \beta_t d\hat W_t^{\Gg} - \frac{1}{2} \int_0^T |\beta_t|^2 dt \Big),
\enqs
under which, by Girsanov's theorem (see \cite{bierut02} Ch.5.2), $\overline W^{\Gg}$ $=$ $\hat W^{\Gg} + \int \beta dt$ is a $(\Qq,\Gg)$-Brownian motion,  
$M$ is a $(\Qq,\Gg)$-martingale, so that the dynamics of $S$ follows a $(\Qq,\Gg)$-local martingale:
\beqs
dS_t &=& S_{t^-} (\sigma_t d\overline W_t^{\Gg} - \gamma_t dM_t). 
\enqs
We thus have the ``no-arbitrage" condition
\beq \label{probamar}
\Mc(\Gg) &:=& \{ \Qq \sim  \Pp \mbox{ on } (\Omega,\Gc_T): S \mbox{ is a } (\Qq,\Gg)-\mbox{local martingale} \} \; \neq \; \emptyset. 
\enq
}
\end{Rem}

\section{Decomposition of  the  utility maximization problem}

\setcounter{equation}{0}
\setcounter{Thm}{0}
\setcounter{Pro}{0}
\setcounter{Cor}{0}
\setcounter{Lem}{0}
\setcounter{Def}{0}
\setcounter{Rem}{0}
\setcounter{Exe}{0}

We are given an utility function $U$ defined on $(0,\infty)$, strictly increasing, strictly concave and $C^1$ on $(0,\infty)$, and satisfying the Inada conditions
$U'(0^+)$ $=$ $\infty$, $U'(\infty)$ $=$ $0$. The performance of an admissible trading strategy
$\pi$ $\in$ $\Ac$ associated to a wealth process $X$ solution to \reff{dynX} and starting at time $0$ from $X_0$ $>$  $0$,  is measured over the finite horizon $T$ by:
\beqs
J_0(\pi) &=& \Ee [U(X_T)],
\enqs
and the optimal investment problem is formulated as:
\beq \label{defv}
V_0  &=& \sup_{\pi\in\Ac} J_0(\pi).
\enq
Problem \reff{defv} is a maximization problem of expected utility from terminal wealth in an  incomplete market  due to the jump of the risky asset.
This optimization problem can be studied by convex duality methods.  Actually, under the condition that
\beq \label{V0finite}
V_0 & < & \infty,
\enq
which is satisfied  under \reff{probamar} once
\beqs
\Ee\Big[ \tilde U\big(y \frac{d\Qq}{d\Pp} \big) \Big] & < & \infty, \;\;\; \mbox{ for some } y >0,
\enqs
where $\tilde U(y)$ $=$ $\sup_{x>0}[U(x)-xy]$, and  under  the  so-called condition of  reasonable asymptotic elasticity:
\beqs
AE(U) := \limsup_{x\rightarrow\infty} \frac{xU'(x)}{U(x)} \; < \; 1,
\enqs
we know from the general results of  Kramkov and Schachermayer \cite{krasch99} that  there exists a solution to \reff{defv}.
We also have a dual characterization of the solution, but this does
not lead to explicit results due to the incompleteness of the market, i.e. the infinite cardinality of $\Mc(\Gg)$.  One can also deal with problem  \reff{defv}
by dynamic programming methods as done recently in Lim and Quenez \cite{Lim-Quenez} under {\bf (H)} hypothesis, 
but again, except for the logarithmic utility function, this does not yield explicit characterization of the optimal strategy. 
We provide here an alternative approach  by  making use of the specific feature of the jump of the stock induced by the default time under the density hypothesis.   The main idea is to separate the problem in two portfolio optimization  problems in complete markets: the after-default and before-default
maximization problems.  This gives a better understanding of the optimal strategy and  
allows us to derive  explicit results in some particular cases of interest.

The  derivation starts as follows.  First notice that any $\pi$ $\in$ $\Ac$, thus  in the form:
$\pi_t$ $=$ $\pi_t^{\Ff}1_{t\leq\tau}$ $+$ $\pi_t^d(\tau)1_{t>\tau}$, $0\leq t\leq T$, can be identified with a pair $(\pi^{\Ff},\pi^d)$ $\in$
$\Ac_{\Ff}\times\Ac_d$ where $\Ac_{\Ff}$ is  the set of  admissible trading strategies   in absence of  defaults, i.e. the set of $\Ff$-adapted processes
$\pi^{\Ff}$ s.t.
\beq \label{defAF}
\int_0^T |\pi_t^{\Ff} \sigma_t^{\Ff}|^2 dt  \; < \;  \infty,  & \mbox{ and } & \pi_\theta \gamma_\theta <  1,  \;\;\; 0 \leq \theta \leq T, \; a.s.
\enq
and $\Ac_d$ is the set of admissible trading strategies after default at time $\tau$ $=$ $\theta$, i.e. the set of Borel family of
$\Ff$-adapted processes
$\{\pi^d_t(\theta), \theta< t \leq T\}$ parametrized by $\theta$ $\in$ $[0,T]$ s.t.
\beqs
\int_\theta^T |\pi_t^d(\theta)\sigma_t^d(\theta)|^2 dt & < & \infty, \;\;\;\;\; a.s.
\enqs
Hence,  for any $\pi$ $=$ $(\pi^{\Ff},\pi^d)$ $\in$ $\Ac$, we observe  by  \reff{Xbefore}-\reff{Xat}-\reff{Xafter}
that the  terminal wealth is written as:
\beqs
X_T &=& X_T^{\Ff} 1_{\tau > T} + X_T^d(\tau)  1_{\tau\leq T},
\enqs
where $X^{\Ff}$ is the wealth process in absence of default, governed by:
\beq \label{dynXF}
dX_t^{\Ff} &=& X_t^{\Ff} \pi_t^{\Ff}   \big( \mu_t^{\Ff} dt + \sigma_t^{\Ff} dW_t\big), \;\;\;\;\; 0 \leq t \leq T,
\enq
starting from $X_0^{\Ff}$ $=$ $X_0$, and
$\{X_t^d(\theta), \theta\leq t \leq T\}$  is the wealth process after default occuring  at $\tau$ $=$ $\theta \in [0,T]$, governed by:
\beq
dX_t^d(\theta) &=& X_t^d(\theta) \pi_t^d(\theta)  \big(  \mu_t^d(\theta) dt + \sigma_t^d(\theta) dW_t \big), \;\;\;\;\; \theta < t \leq T \label{dynXa} \\
X_\theta^d(\theta) & = &  X_\theta^{\Ff}(1-\pi_\theta^{\Ff} \gamma_\theta).
\enq
Therefore, under the density hypothesis {\bf (DH)}, and by the law of iterated conditional expectations,  the performance measure may be written as:
\beq
J_0(\pi) &=& \Ee \big[ \Ee [U(X_T)| \Fc_T] \big] \; = \; \Ee\big[  U(X_T^{\Ff}) \Pp[\tau > T | \Fc_T] +  \Ee[U(X_T^d(\tau)) 1_{\tau \leq T} | \Fc_T] \big] \nonumber \\
&=& \Ee \Big[  U(X_T^{\Ff}) G_T  +  \int_0^T U(X_T^d(\theta))  \alpha_T(\theta) d\theta \Big], \label{Jdec}
\enq
where  $G_T$ $=$ $\Pp[\tau > T | \Fc_T]$ $=$ $\int_T^\infty \alpha_T(\theta) d\theta$.

Let us  introduce  the value-function  process  of the  ``after-default" optimization problem:
\beq
V^d_\theta(x) &=& \esssup_{\pi^d(\theta)\in\Ac_{d}(\theta)} J^d_\theta(x,\pi^d(\theta)) , \;\;(\theta,x) \in [0,T]\times (0,\infty), \label{defva} \\
J^d_\theta(x,\pi^d(\theta))  &=& \Ee \big[ U(X_T^{d,x}(\theta)) \alpha_T(\theta) \big| \Fc_\theta \big], \nonumber
\enq
where $\Ac_{d}(\theta)$ is the  set of  $(\Fc_t)_{\theta<t\leq T}$-adapted processes  $\{\pi^d_t(\theta), \theta< t \leq T\}$ satisfying
$\int_\theta^T |\pi_t^d(\theta)\sigma_t^d(\theta)|^2dt$ $<$ $\infty$ a.s., and  $\{X_t^{d,x}(\theta), \theta\leq t\leq T\}$ is the solution to \reff{dynXa} controlled by  $\pi^d(\theta)$ $\in$ $\Ac_d(\theta)$, starting from $x$ at time $\theta$.  Thus,
$V^d$ is the value-function process  of an optimal investment problem in a market model after default.  Notice that the coefficients $(\mu^d,\sigma^d)$
of the model depend on the initial time $\theta$ when the maximization is performed, and the utility function in the criterion is weighted by
$\alpha_T(\theta)$.  We shall see in the next section how to deal with these peculiarities  for solving \reff{defva} and  proving the existence and characterization of  an optimal strategy.

The main result of this section is to show that the original problem \reff{defv}  can be split into the above after-default optimization problem, and a global optimization problem in a before-default market.

\begin{Thm}
Assume that  $V_\theta^d(x)$ $<$ $\infty$ a.s. for all $(\theta,x)\in [0,T]\times (0,\infty)$. Then,  we have:
\beq \label{vglobal}
V_0 &=& \sup_{\pi^{\Ff}\in\Ac_{\Ff}} \Ee \Big[ U(X_T^{\Ff}) G_T
+  \int_0^T V^d_\theta(X_\theta^{\Ff}(1- \pi_\theta^{\Ff} \gamma_\theta)) d\theta \Big].
\enq
\end{Thm}
{\bf Proof.} Given $\pi$ $=$ $(\pi^{\Ff},\pi^d)$ $\in$ $\Ac$,  we have the relation \reff{Jdec} for $J_0(\pi)$ under {\bf (DH)}.  Furthermore,  by Fubini's theorem,  the law of iterated conditional expectations, we then obtain:
\beq
J_0(\pi) &=&
\Ee \Big[  U(X_T^{\Ff}) G_T  +  \int_0^T \Ee \big[  U(X_T^d(\theta))  \alpha_T(\theta) \big| \Fc_\theta]  d\theta \Big]  \nonumber \\
&=& \Ee \Big[  U(X_T^{\Ff}) G_T  + \int_0^T J^d_\theta(X_\theta^d(\theta),\pi^d(\theta)) d \theta \Big] \label{Jd}  \\
& \leq &  \Ee \Big[  U(X_T^{\Ff}) G_T  + \int_0^T V^d_\theta(X_\theta^d(\theta)) d \theta \Big]  \nonumber \\
& \leq & \sup_{\pi^{\Ff}\in\Ac_{\Ff}} \Ee \Big[ U(X_T^{\Ff}) G_T
+  \int_0^T V^d_\theta(X_\theta^{\Ff}(1- \pi_\theta^{\Ff} \gamma_\theta)) d\theta \Big] =: \hat V_0. \nonumber
\enq
by definitions of $J^d$, $V^d$ and $X_\theta^d(\theta)$.
This proves the inequality: $V_0$ $\leq$ $\hat V_0$.

To prove the converse inequality, fix an arbitrary $\pi^{\Ff}$ $\in$ $\Ac_{\Ff}$.
By definition of   $V^d$, for any $\omega$ $\in$ $\Omega$, $\theta$ $\in$ $[0,T]$,  and $\eps$ $>$ $0$,
there exists $\pi^{d,\eps,\omega}(\theta)$ $\in$ $\Ac_d(\theta)$, which is an $\eps$-optimal control
for $V^d_\theta$ at  $(\omega,X_\theta^d(\omega,\theta))$.  By a measurable selection result (see e.g. \cite{wag80}), 
one can find $\pi^{d,\eps}$ $\in$ $\Ac_d$ s.t.  $\pi^{d,\eps}(\omega,\theta)$ $=$ $\pi^{d,\eps,\omega}(\omega,\theta)$, 
$d\Pp\otimes d\theta$ a.e., and so 
\beqs
V^d_\theta(X_\theta^d(\theta)) - \eps & \leq &  J^d_\theta(X_\theta^d(\theta),\pi^{d,\eps}(\theta)), \;\;\; d\Pp\otimes d\theta \; a.e. 
\enqs
By denoting $\pi^\eps$ $=$ $(\pi^{\Ff},\pi^{d,\eps})$ $\in$ $\Ac$, and using again  \reff{Jd}, we then get:
\beqs
V_0  & \geq & J_0(\pi^\eps) \; = \;
\Ee \Big[  U(X_T^{\Ff}) G_T  + \int_0^T J^d_\theta(X_\theta^d(\theta),\pi^{d,\eps}(\theta)) d \theta \Big]   \\
& \geq & \Ee \Big[  U(X_T^{\Ff}) G_T  + \int_0^T V^d_\theta(X_\theta^d(\theta)) d \theta \Big] - \eps.
\enqs
From the arbitrariness of $\pi^{\Ff}$ in $\Ac_{\Ff}$ and $\eps$ $>$ $0$, we obtain the required inequality and so the result.
\ep

\begin{Rem} \label{remdec}
{\rm
The relation \reff{vglobal} can be viewed as a dynamic programming type relation. Indeed, as in dynamic programming principle (DPP),
we look for a relation on  the value function by varying the initial states. However,  instead of taking  two consecutive dates as in the usual
DPP, the original feature here is to derive the equation by considering
the value function between the initial time and the default time conditionnally on  the terminal information,  leading to the introduction of an  ``after-default"
and a global before-default optimization problem, the latter  involving the former.  Each of these optimization problems are performed in complete market models driven by the brownian motion and with coefficients adapted with respect to the brownian filtration.
The main advantage of this approach is then to reduce the problem to the resolution of two optimization problems in complete markets,
which are simpler to  deal with, and give more explicit results  than the incomplete market  framework studied by the   ``classical" dynamic programming approach or the convex duality method.

Furthermore,  a careful  look at the arguments for deriving the relation \reff{vglobal} shows that in the decomposition of  the optimal trading strategy for the original problem \reff{defv} which is known to exist a priori under \reff{probamar}:
\beqs
\hat \pi_t &=& \hat \pi_t^{\Ff} 1_{t\leq \tau} + \hat \pi_t^d(\tau) 1_{t>\tau}, \;\;\; 0 \leq t\leq T,
\enqs
$\hat\pi^{\Ff}$ is an optimal control to \reff{vglobal}, and $\hat\pi^d(\theta)$ is an optimal control to  $V^d_\theta(\hat X_\theta^d(\theta))$
with  $\hat X_\theta^d(\theta)$ $=$ $\hat X_\theta^{\Ff}(1- \hat \pi_\theta^{\Ff} \gamma_\theta)$, and  $\hat X^{\Ff}$
is the wealth process governed by $\hat\pi^{\Ff}$. In other words, the optimal trading strategy is to follow the trading strategy $\hat\pi^{\Ff}$ before default
time $\tau$, and then to change to the after-default trading strategy $\hat\pi^d(\tau)$, which depends on the time where default occurs.
In the next section, we focus on the resolution of these two optimization problems.
}
\end{Rem}

\section{Solution to the optimal investment problem}

\setcounter{equation}{0}
\setcounter{Thm}{0}

In this section, we focus on the resolution of the two optimization problems arising from the decomposition of  the initial utility maximization problem.
We first study the after-default optimal invesment problem, and then the global before-default optimization problem.

\subsection{The after-default utility maximization problem}\label{Subsec: after default optimization}

Problem \reff{defva} is an optimal investment problem in a complete market model after default.
A specific feature of this model is the dependence of the coefficients
$(\mu^d,\sigma^d)$ on the initial time $\theta$ when the  maximization is performed. This makes the optimization problem time-inconsistent, and the classical dynamic programming method can not be applied.  Another peculiarity in the criterion is the presence of the density term $\alpha_T(\theta)$
weighting the utility function $U$.

We  adapt the  convex duality method  for solving \reff{defva}.  We have to extend this martingale  method (in complete market)
in a dynamic framework, since we want to compute the value-function process  at any time $\theta$ $\in$ $[0,T]$. Let us denote by:
\beqs
Z_t(\theta) &=& \exp\Big(  -\int_\theta^t  \frac{\mu_u^d(\theta)}{\sigma_u^d(\theta)} dW_u
- \frac{1}{2} \int_\theta^t \Big|\frac{\mu_u^d(\theta)}{\sigma_u^d(\theta)}\Big|^2 du \Big), \;\;\;  \theta\leq t\leq T,
\enqs
the (local) martingale density in the market model \reff{Safter} after default.  We assume that  for all $\theta$ $\in$ $[0,T]$, 
 there exists some $y_\theta$ $\Fc_\theta$-measurable strictly positive random variable  s.t.
\beq \label{tildeUZtheta}
\Ee \Big[ \tilde U\Big(y_\theta \frac{Z_T(\theta)}{\alpha_T(\theta)}\Big) \alpha_T(\theta) \Big| \Fc_\theta \Big] & < & \infty.
\enq
This assumption is similar to the one imposed in the classical (static) convex duality method for ensuring that the dual problem is well-defined and finite.

\begin{Thm} \label{thmafter}
Assume that \reff{tildeUZtheta} and $AE(U)$ $<$ $1$ hold true.
Then, the value-function process  to problem \reff{defva} is  finite a.s.  and given by
\beqs
V^d_\theta(x) &=&   \Ee \Big[ U \Big( I\Big( \hat y_\theta(x) \frac{Z_T(\theta)}{\alpha_T(\theta)}  \Big) \Big)\alpha_T(\theta) \Big | \Fc_\theta  \Big], \;\;\;
(\theta,x)  \in [0,T]\times (0,\infty),
\enqs
and the corresponding optimal wealth process is equal to:
\beq \label{defhatX}
\hat X_t^{d,x}(\theta) &=& \Ee \Big[ \frac{Z_T(\theta)}{Z_t(\theta)}  I\Big( \hat y_\theta(x) \frac{Z_T(\theta)}{\alpha_T(\theta)}  \Big) \Big | \Fc_t \Big], \;\;\;
\theta \leq t\leq T,
\enq
where $I$ $=$ $(U')^{-1}$ is the inverse of $U'$, and $\hat y_\theta(x)$ is the strictly positive $\Fc_\theta\otimes\Bc((0,\infty))$-measurable random variable
solution to $\hat X_\theta^{d,x}(\theta)$ $=$ $x$. 
\end{Thm}
{\bf Proof.}   First observe, similarly as in Theorem 2.2  in \cite{krasch99}, that under  $AE(U)$ $<$ $1$, the validity of \reff{tildeUZtheta} for some or for all
$y_\theta$ $\Fc_\theta$-measurable strictly positive random variable, is equivalent.
By definition of $Z(\theta)$ and It\^o's formula, the process $\{Z_t(\theta)X_t^{d,x}(\theta), \theta\leq t\leq T\}$ is a nonnegative
$(\Pp,(\Fc_t)_{\theta\leq t\leq T})$-local martingale, hence a supermartingale,  for any $\pi^d(\theta)$ $\in$ $\Ac_d(\theta)$, and so
$\Ee[ X_T^{d,x}(\theta)  Z_T(\theta) | \Fc_\theta]$ $\leq$ $X_\theta^{d,x}(\theta)Z_\theta(\theta)$ $=$ $x$.  Denote
$Y_T(\theta)$ $=$ $Z_T(\theta)/\alpha_T(\theta)$.  Then,  by definition  of $\tilde U$, we have for all  $y_\theta$  $\Fc_\theta$-measurable strictly positive random variable, and $\pi^d(\theta)$ $\in$  $\Ac_d(\theta)$:
\beq
\Ee[ U(X_T^{d,x}(\theta)) \alpha_T(\theta) | \Fc_\theta]  & \leq &
\Ee[ \tilde U(y_\theta Y_T(\theta)) \alpha_T(\theta) | \Fc_\theta]  +  \Ee[ X_T^{d,x}(\theta) y_\theta Y_T(\theta) \alpha_T(\theta) | \Fc_\theta]
\nonumber \\
&=&   \Ee[ \tilde U(y_\theta Y_T(\theta)) \alpha_T(\theta) | \Fc_\theta]  +  y_\theta   \Ee[ X_T^{d,x}(\theta)  Z_T(\theta) |\Fc_\theta ] \nonumber \\
& \leq &  \Ee[ \tilde U(y_\theta Y_T(\theta)) \alpha_T(\theta) | \Fc_\theta ]  + xy_\theta,  \label{interVthetafini}
\enq
which proves in particular that $V_\theta^d(x)$ is finite a.s.  Now,  we recall that under the Inada conditions,
the supremum in the definition of  $\tilde U(y)$  is attained at $I$, i.e. $\tilde U(y)$ $=$ $U(I(y))-yI(y)$. From \reff{interVthetafini}, this implies
\beq
\Ee[ U(X_T^{d,x}(\theta)) \alpha_T(\theta) | \Fc_\theta]  & \leq &
\Ee[ U(I(y_\theta Y_T(\theta))) \alpha_T(\theta)| \Fc_\theta]  \nonumber \\
& & \;\; - \;  y_\theta \Big(\Ee[Z_T(\theta)   I(y_\theta Y_T(\theta)) |\Fc_\theta ] - x \Big).  \label{interdual}
\enq
Now, under the Inada conditions, \reff{tildeUZtheta} and $AE(U)$ $<$ $1$, for any $\omega$ $\in$ $\Omega$, $\theta$ $\in$ $[0,T]$,
the function $y$ $\in$ $(0,\infty)$ $\rightarrow$ $f_\theta(\omega,y)$ $=$
$\Ee[Z_T(\theta)   I(yY_T(\theta)) |\Fc_\theta]$ is a strictly decreasing
one-to-one continuous function from $(0,\infty)$ into $(0,\infty)$.  Hence, there exists a unique $\hat y_\theta(\omega,x)$ $>$ $0$ s.t.
$f_\theta(\omega,\hat y(x))$ $=$ $x$.  Moreover, since $f_\theta(y)$ is $\Fc_\theta\otimes\Bc(0,\infty)$-measurable, this value $\hat y_\theta(x)$ can be chosen,  by a measurable selection argument, as $\Fc_\theta\otimes\Bc(0,\infty)$-measurable. With this choice of $y_\theta$ $=$ $\hat y_\theta(x)$,
and by setting $\hat X_T^{d,x}(\theta)$ $=$ $I(\hat y_\theta(x) Y_T(\theta))$,  the inequality  \reff{interdual} yields:
\beq \label{Uless}
\Ee[ U(X_T^{d,x}(\theta)) \alpha_T(\theta) |\Fc_\theta]  & \leq &
\Ee[ U(\hat X_T^{d,x}(\theta)) \alpha_T(\theta) |\Fc_\theta], \;\;\; \forall \pi^d(\theta) \in \Ac_d(\theta).
\enq
Consider now the process $\hat X^{d,x}(\theta)$ defined  in \reff{defhatX} leading to $\hat X_T^{d,x}(\theta)$ at time $T$.
By definition, the process
$\{M_t(\theta)=Z_t(\theta)\hat X_t^{d,x}(\theta), \theta\leq t\leq T\}$ is a strictly positive $(\Pp,(\Fc_t)_{\theta\leq t\leq T})$-martingale.
From the martingale representation theorem for brownian motion filtration, there exists an $(\Fc_t)_{\theta\leq t\leq T}$-adapted process
$(\phi_t)_{\theta\leq t\leq T}$ satisfying  $\int_\theta^T |\phi_t|^2 dt$ $<$ $\infty$ a.s., and such that
\beqs
M_t(\theta) &=& M_\theta(\theta) + \int_\theta^t  \phi_u M_u(\theta) dW_u, \;\;\;  \theta \leq t \leq T.
\enqs
Thus, by setting $\hat\pi^d(\theta)$ $=$ $(\phi + \frac{\mu^d(\theta)}{\sigma^d(\theta)})/\sigma^d(\theta)$, we see that $\hat\pi^d(\theta)$ $\in$
$\Ac_d(\theta)$, and by It\^o's formula,  $\hat X^{d,x}(\theta)$ $=$ $M(\theta)/Z(\theta)$  satisfies the  wealth equation \reff{dynXa} controlled by
$\hat\pi^d(\theta)$.  Moreover, by construction of $\hat y_\theta(x)$, we have:
\beqs
\hat X_\theta^{d,x}(\theta) &=&  \Ee \Big[ Z_T(\theta)   I\Big( \hat y_\theta(x) \frac{Z_T(\theta)}{\alpha_T(\theta)}  \Big) \Big| \Fc_\theta  \Big]  \; = \;  x.
\enqs
Recalling  \reff{Uless}, this proves that $\hat\pi^d(\theta)$ is an optimal solution to   \reff{defva}, with corresponding
optimal wealth process $\hat X^{d,x}(\theta)$.
\ep

\begin{Rem}
{\rm Under the {\bf (H)} hypothesis, $\alpha_T(\theta)$ $=$ $\alpha_\theta(\theta)$ is $\Fc_\theta$-measurable. In this case, the optimal wealth process to
\reff{defva} is given by:
\beqs
\hat X_t^{d,x}(\theta) &=& \Ee \Big[ \frac{Z_T(\theta)}{Z_t(\theta)}  I\big( \bar y_\theta(x) Z_T(\theta) \big) \Big | \Fc_t \Big], \;\;\;
\theta \leq t\leq T,
\enqs
where $\bar y_\theta(x)$ is the strictly positive $\Fc_\theta\otimes\Bc((0,\infty))$-measurable random variable satisfying $\hat X_\theta^{d,x}(\theta)$ $=$
$x$.  Hence, the optimal strategy after-default does not depend on the density of the default time.
}
\end{Rem}

\vspace{2mm}

We illustrate the above results in the case of Constant Relative Risk Aversion (CRRA) utility functions.

\begin{Exe}{\rm {\bf The case of CRRA Utility function}

\vspace{1mm}

\noindent We consider  utility functions  in the form
\beqs
U(x) &=& \frac{x^p}{p}, \;\;\;\;\;  p<1, p \neq 0, \; x >0.
\enqs
In this case, we easily compute the optimal wealth process in \reff{defhatX}:
\beqs
\hat X_t^{d,x}(\theta) &=& \frac{x}{\Ee\Big[\alpha_T(\theta)  \Big(\frac{Z_T(\theta)}{\alpha_T(\theta)} \Big)^{-q} \Big| \Fc_\theta\Big]}.
\frac{ \Ee\Big[\alpha_T(\theta)   \Big(\frac{Z_T(\theta)}{\alpha_T(\theta)} \Big)^{-q}  \Big| \Fc_t\Big]}{Z_t(\theta)}, \;\;\; \theta\leq t \leq T,
\enqs
where  $q$ $=$ $\frac{p}{1-p}$. The optimal value process is then given for all $x$ $>$ $0$ by:
\beq \label{Vdpower}
V^d_\theta(x) &=& \frac{x^p}{p}.  \Big( \Ee\Big[\alpha_T(\theta)  \Big(\frac{Z_T(\theta)}{\alpha_T(\theta)} \Big)^{-q} \Big| \Fc_\theta\Big] \Big)^{1-p},
\;\;\; \theta \in [0,T].
\enq
Notice that the case of logarithmic utility function:  $U(x)$ $=$ $\ln x$, $x$ $>$ $0$, can be either computed directly, or derived as the limiting case
of power utility function case:  $U(x)$ $=$ $\frac{x^p-1}{p}$ as $p$ goes to zero. The optimal wealth process is given by:
\beqs
\hat X_t^{d,x}(\theta) &=& \frac{x}{\Ee[\alpha_T(\theta) | \Fc_\theta]}. \frac{ \Ee[\alpha_T(\theta) | \Fc_t]}{Z_t(\theta)}, \;\;\; \theta\leq t \leq T,
\enqs
and the optimal value process for all $x$ $>$ $0$, is equal to:
\beqs
V^d_\theta(x) &=& \Ee[\alpha_T(\theta) | \Fc_\theta] \ln \Big(\frac{x}{\Ee[\alpha_T(\theta) | \Fc_\theta]} \Big) +
\Ee\Big[\alpha_T(\theta)  \ln \Big( \frac{\alpha_T(\theta)}{Z_T(\theta)}\Big) \Big| \Fc_\theta\Big], \;\;\; \theta \in [0,T].
\enqs
}
\end{Exe}

\subsection{The global before-default optimization problem}

In this paragraph, we focus on the resolution of the optimization problem \reff{vglobal}. We already know the existence of an optimal strategy
$\hat\pi^{\Ff}$ to  this problem, see Remark \ref{remdec}, and our main concern is to provide an explicit characterization of the optimal control.

We use a dynamic programming approach.  For any $t$ $\in$ $[0,T]$, $\nu$ $\in$ $\Ac_{\Ff}$,  let us consider  the set of controls coinciding with $\nu$
until time $t$:
\beqs
\Ac_{\Ff}(t,\nu) &=& \{ \pi^{\Ff} \in \Ac_{\Ff}: \pi^{\Ff}_{.\wedge t} = \nu_{.\wedge t} \}.
\enqs
Under the standing condition that $V_0$ $<$ $\infty$, 
we then introduce  the dynamic version of the optimization problem \reff{vglobal} by  considering the family of
$\Ff$-adapted processes:
\beqs
V_t(\nu) &=& \esssup_{\pi^{\Ff} \in \Ac_{\Ff}(t,\nu)} \Ee  \Big[ U(X_T^{\Ff}) G_T
+  \int_t^T V^d_\theta(X_\theta^{\Ff}(1- \pi_\theta^{\Ff} \gamma_\theta))  d\theta \Big| \Fc_t \Big], \;\;\; 0 \leq t\leq T,
\enqs
so that $V_0$ $=$ $V_0(\nu)$ for any $\nu$ $\in$ $\Ac_{\Ff}$.  In the above expression, $X^{\Ff}$ is the wealth process of dynamics  \reff{dynXF},
controlled by $\pi^{\Ff}$ $\in$ $\Ac(t,\nu)$,  and starting from $X_0$.  We also denote $X^{\nu,\Ff}$ the wealth process of dynamics  \reff{dynXF},
controlled by $\nu$ $\in$ $\Ac_{\Ff}$, starting from $X_0$, so that it coincides with $X^{\Ff}$ until time $t$, i.e.
$X_{.\wedge t}^{\nu,\Ff}$ $=$ $X_{. \wedge t}^{\Ff}$.
From the dynamic programming principle (see El Karoui \cite{elk81}),
the process $\{V_t(\nu),0\leq t\leq T\}$ can be chosen in its c\`ad-l\`ag version, and is such that  for any $\nu$ $\in$ $\Ac_{\Ff}$:
\beq \label{dynprosurmar}
\big\{V_t(\nu) + \int_0^t  V^d_\theta(X_\theta^{\nu,\Ff}(1- \nu_\theta \gamma_\theta))  d\theta, \; 0\leq t\leq T\big\}
& \mbox{is a} & (\Pp,\Ff)-\mbox{supermartingale}.
\enq
Moreover,  the optimal strategy $\hat\pi^{\Ff}$ to problem $V_0$, is characterized by the martingale pro\-perty:
\beq \label{dynpromar}
\big\{V_t(\hat \pi^{\Ff}) + \int_0^t  V^d_\theta(X_\theta^{\hat\pi^{\Ff},\Ff}(1- \hat\pi_\theta^{\Ff} \gamma_\theta))  d\theta, \;  0 \leq t\leq T \big\}
& \mbox{is a} & (\Pp,\Ff)-\mbox{martingale}.
\enq

In the sequel, we shall exploit this dynamic programming properties in the  particular  important case of constant relative risk aversion (CRRA)  
utility functions. 
We then consider  utility functions in the form 
\beqs
U(x) &=& \frac{x^p}{p}, \;\;\;\;\;  p<1, \; p \neq 0,   \; x >0,
\enqs
and we set $q$ $=$ $\frac{p}{1-p}$.  Notice that we deal with the relevant economic case when  $p$ $<$ $0$, i.e.  
the degree of risk aversion $1-p$ is strictly larger than $1$.  This  will induce  some additional technical difficulties with respect  to the case 
$p$ $>$ $0$. 
For CRRA utility function,  $V^d(x)$ is also of the same power type, see \reff{Vdpower}:
\beqs
V_\theta^d(x) &=&  U(x)   K_\theta^p \;\;\; \mbox{ with }  \;\;\;
K_\theta \; = \;  \Big( \Ee\Big[\alpha_T(\theta)  \Big(\frac{Z_T(\theta)}{\alpha_T(\theta)} \Big)^{-q} \Big| \Fc_\theta\Big] \Big)^{1\over q},
\enqs
and we assume that $K_\theta$ is   finite a.s. for all $\theta$ $\in$ $[0,T]$. The value of the optimization problem \reff{vglobal} is written as
\beqs
V_0 &=& \sup_{\nu\in\Ac_{\Ff}} \Ee[ U(X_T^{\nu,\Ff}) G_T + \int_0^T U(X_\theta^{\nu,\Ff})(1-\nu_\theta\gamma_\theta)^p K_\theta^p d\theta],
\enqs
In the above equality, we may without loss of generality take supremum over $\Ac_{\Ff}(U)$, the set of elements $\nu$ $\in$ $\Ac_{\Ff}$ such that: 
\beq \label{integV0}
\Ee[ U(X_T^{\nu,\Ff}) G_T + \int_0^T U(X_\theta^{\nu,\Ff})(1-\nu_\theta\gamma_\theta)^p K_\theta^p d\theta] & > &- \infty,
\enq
and by misuse of notation, we write  $\Ac_{\Ff}$ $=$ $\Ac_{\Ff}(U)$. For any $\nu$ $\in$ $\Ac_{\Ff}$ with corresponding
wealth process $X^{\nu,\Ff}$ governed by \reff{dynXF} with control $\nu$, and starting from $X_0$,  we notice that the c\`ad-l\`ag
$\Ff$-adapted
process defined by:
\beq
Y_t &:= & \frac{V_t(\nu)}{U(X_t^{\nu,\Ff})} \label{defYdynpro} \\
& = &  p \esssup_{\pi^{\Ff} \in \Ac_{\Ff}(t,\nu)} \Ee  \Big[ U\Big(\frac{X_T^{\Ff}}{X_t^{\nu,\Ff}}\Big) G_T
+  \int_t^T U\Big(\frac{X_\theta^{\Ff}}{X_t^{\nu,\Ff}}\Big)(1- \pi_\theta^{\Ff} \gamma_\theta)^p K_\theta^p  d\theta \Big| \Fc_t \Big], \;\;\;  0 \leq t\leq T \nonumber
\enq
does not depend on $\nu$ $\in$ $\Ac_{\Ff}$.  It lies in the set $L_+(\Ff)$ of nonnegative c\`ad-l\`ag $\Ff$-adapted processes.  Let us also denote by
$L_{loc}^2(W)$ the set of $\Ff$-adapted process $\phi$ s.t. $\int_0^T |\phi_t|^2 dt$ $<$ $\infty$ a.s.

\vspace{1mm}

We have the following preliminary properties  on this  process $Y$.

\begin{Lem} \label{lempre} 
The process $Y$ in \reff{defYdynpro} is strictly positive, i.e.  $\Pp[Y_t > 0, \; 0\leq t\leq T]$ $=$ $1$. 
Moreover, for all $\nu$ $\in$ $\Ac_{\Ff}$, the process 
\beq \label{defxinuY}
\xi_t^\nu(Y) &:=&  U(X_t^{\nu,\Ff}) Y_t + \int_0^t U(X_\theta^{\nu,\Ff}) (1-\nu_\theta\gamma_\theta)^p K_\theta^p d\theta, \;\;\; 0 \leq t\leq T, 
\enq
is  bounded from below by a martingale. 
\end{Lem}
{\bf Proof.} (1) We first consider the case $p$ $>$ $0$. Then, 
\beq
Y_t &= &  \esssup_{\pi^{\Ff} \in \Ac_{\Ff}(t,\nu)} \Ee  \Big[ \Big(\frac{X_T^{\Ff}}{X_t^{\nu,\Ff}}\Big)^p G_T
+  \int_t^T \Big(\frac{X_\theta^{\Ff}}{X_t^{\nu,\Ff}}\Big)^p(1- \pi_\theta^{\Ff} \gamma_\theta)^p K_\theta^p  d\theta \Big| \Fc_t \Big] \label{Yppos} \\
& \geq &  \Ee \big[  G_T + \int_t^T   K_\theta^p d\theta \Big|\Fc_t \big]  \; > \; 0, \;\;\;\; \forall \; t \in [0,T], \nonumber
\enq
by taking in \reff{Yppos}  the control process $\pi^{\Ff}$ $\in$ $\Ac_{\Ff}(t,\nu)$ defined by $\hat\pi^{\Ff}_s$ $=$ $\nu_s 1_{s\leq t}$. Moreover, since 
$U(x)$ is nonnegative, the process $\xi^\nu(Y)$ is nonnegative, hence trivially bounded from below by the zero  martingale. 

\noindent (2) We next consider the case $p$ $<$ $0$. Then, 
\beq
Y_t &= &  \essinf_{\pi^{\Ff} \in \Ac_{\Ff}(t,\nu)} \Ee  \Big[ \Big(\frac{X_T^{\Ff}}{X_t^{\nu,\Ff}}\Big)^p G_T
+  \int_t^T \Big(\frac{X_\theta^{\Ff}}{X_t^{\nu,\Ff}}\Big)^p(1- \pi_\theta^{\Ff} \gamma_\theta)^p K_\theta^p  d\theta \Big| \Fc_t \Big]   \label{Ypneginter} \\
& \geq &   J_t :=  \essinf_{\pi^{\Ff} \in \Ac_{\Ff}(t,\nu)} \Ee  \Big[ \Big(\frac{X_T^{\Ff}}{X_t^{\nu,\Ff}}\Big)^p G_T \Big| \Fc_t \Big], \;\;\; \forall \; t \in [0,T].  \nonumber 
\enq
Notice that the process $J$ can be chosen in its c\`ad-l\`ag modification. 
Let us show that for any $t$ $\in$ $[0,T]$,  the infimum in $J_t$ is attained.  
Fix $t$ $\in$ $[0,T]$, and 
consider, by a measurable selection argument,  a minimizing sequence $(\pi^n)_n$ $\in$ $\Ac_{\Ff}(t,\nu)$ to $J_t$, i.e.
\beq \label{limJ}
\lim_{n\rightarrow\infty} \Ee  \Big[ \Big(\frac{X_T^{n}}{X_t^{\nu,\Ff}}\Big)^p G_T \Big| \Fc_t \Big] &=& J_t, \;\;\; a.s. 
\enq
Here $X^{n}$ denotes the wealth process of dynamics \reff{dynXF} governed by $\pi^n$.  Consider the (local) martingale density process 
\beqs
Z_s^t &=& \exp\Big( - \int_t^s \frac{\mu_u^{\Ff}}{\sigma_u^{\Ff}} dW_u - \frac{1}{2} \int_t^s \Big| \frac{\mu_u^{\Ff}}{\sigma_u^{\Ff}} \Big|^2 du \Big), 
\;\;\; t \leq s \leq T. 
\enqs
By definition of $Z^t$ and It\^o's formula, the process $\{ Z_s^t X_s^{n}, t\leq s \leq T\}$ is a nonnegative $(\Pp,(\Fc_s)_{t\leq s\leq T})$-local martingale, hence a supermartingale, and so $\Ee[X_T^{n}Z_T^t | \Fc_t]$ $\leq$ $X_t^{n}Z_t^t$ $=$ $X_t^{\nu,\Ff}$. 
By Koml\`os Lemma  applied to the sequence of nonnegative 
$\Fc_T$-measurable random variable $(X_T^{n})_n$, there exists a convex combination $\tilde X_T^n$ $\in$ conv$(X_T^{k,\Ff},k\geq n\}$  such that 
$(\tilde X_T^n)_n$ converges a.s. to some nonnegative   $\Fc_T$-measurable random variable $\tilde X_T$.  
By Fatou's lemma, we  have $\tilde X_t$ $:=$ $\Ee[\tilde X_T Z_T^t|\Fc_t]$ $\leq$ $X_t^{\nu,\Ff}$.  Moreover, by convexity of $x$ $\rightarrow$ $x^p$,  
and Fatou's lemma, it follows from \reff{limJ} that 
\beq \label{Jeq}
J_t & \geq &   \Ee  \Big[ \Big(\frac{\tilde X_T}{X_t^{\nu,\Ff}}\Big)^p G_T \Big| \Fc_t \Big], \;\;\; a.s.
\enq
Now,  since $p$ $<$ $0$, $J_t$ $<$ $\infty$, and  $G_T$ $>$ $0$ a.s.,  we deduce that  
$\tilde X_T$ $>$ $0$,   and so $\tilde X_t$ $>$ $0$ a.s.  
Consider the process $\bar X_s^t$ $=$ $\frac{X_t^{\nu,\Ff}}{\tilde X_t}\Ee[ \frac{Z_T^t}{Z_s^t} \tilde X_T | \Fc_s]$, $t\leq s\leq T$.  Then, 
$\{Z_s^t \bar X_s^t, t\leq s\leq T\}$ is a strictly positive $(\Pp,(\Fc_s)_{t\leq s\leq T})$-martingale, and by the martingale representation theorem 
for brownian filtration,  using same arguments as in the end of proof of Theorem \ref{thmafter}, we obtain the existence of an 
$(\Fc_s)_{t\leq s\leq T}$-adapted process $\bar\pi^t$ $=$ $(\bar\pi_s^t)_{t \leq s\leq T}$ satisfying $\int_t^T |\bar\pi_s^t\sigma_s^{\Ff}|^2 ds$ $<$ $\infty$, 
such that $\bar X^t$ satisfies the wealth process dynamics \reff{dynXF} with portfolio $\bar\pi^t$ on $(t,T)$, and starting from $\bar X_t$ $=$ $X_t^{\nu,\Ff}$. 
By considering the portfolio process $\bar\pi$ $\in$ $\Ac_{\Ff}(t,\nu)$ defined by $\bar\pi_s$ $=$ $\nu_s1_{s\leq t} + \bar\pi_s^t 1_{s>t}$, for $0\leq s\leq T$,  
and denoting by $X^{\bar\pi,\Ff}$ the corresponding wealth process, it  follows that $X_s^{\bar\pi,\Ff}$ $=$ $\bar X_s^t$ for $t\leq s\leq T$,  and 
in particular $X_T^{\bar\pi,\Ff}$ $=$ $\bar X_T^t$ $=$ $\frac{X_t^{\nu,\Ff}}{\tilde X_t} \tilde X_T$ $\geq$ $\tilde X_T$ a.s.  From \reff{Jeq}, 
the nonincreasing property of $x$ $\rightarrow$ $x^p$, and  definition of $J_t$, we deduce that
\beq \label{tildeJ}
J_t &=&  \tilde J_t := \Ee  \Big[ \Big(\frac{X_T^{\bar\pi,\Ff}}{X_t^{\nu,\Ff}}\Big)^p G_T \Big| \Fc_t \Big], \;\;\; a.s.
\enq
and as a byproduct that $X_T^{\bar\pi,\Ff}$ $=$ $\tilde X_T$. The  equality \reff{tildeJ} 
means that  the process $J$ $=$ $(J_t)_{t\in [0,T]}$ is a modification of the process $\tilde J$ $=$ $(\tilde J_t)_{t\in [0,T]}$. Since, $J$ and $\tilde J$ are  
c\`ad-l\`ag, they are then indistinguishable, i.e. $\Pp[ J_t = \tilde J_t, \; 0 \leq t\leq T]$ $=$ $1$.  We deduce that the process $J$, and consequently $Y$, inherit  the strict positivity of the process $\tilde J$.

From \reff{Ypneginter}, we have for all $\nu$ $\in$ $\Ac_{\Ff}$, $t$ $\in$ $[0,T]$, 
\beq
\xi_t^\nu(Y) &=&  \esssup_{\pi^{\Ff} \in \Ac_{\Ff}(t,\nu)} \Ee  \Big[  U(X_T^{\Ff}) G_T
+  \int_0^T  U(X_\theta^{\Ff}) (1- \pi_\theta^{\Ff} \gamma_\theta)^p K_\theta^p  d\theta \Big| \Fc_t \Big]  \label{Ypnegxi} \\
& \geq & M_t^\nu :=  \Ee \Big[ U(X_T^{\nu,\Ff}) G_T +  \int_0^T  U(X_\theta^{\nu,\Ff}) (1- \nu_\theta \gamma_\theta)^p K_\theta^p  d\theta \Big| \Fc_t \Big],  
\;\;\;\; t \in [0,T], \nonumber
\enq 
by taking in \reff{Ypnegxi}  the control process $\pi^{\Ff}$ $=$ $\nu$  $\in$ $\Ac_{\Ff}(t,\nu)$.  The negative process $(M_t^\nu)_{t\in [0,T]}$ is an  
integrable (recall \reff{integV0}) martingale,  and the assertions of the  Lemma are proved.  
\ep

\vspace{2mm}

In the sequel, we denote by $L_+^{b}(\Ff)$ the set of processes $\tilde Y$ in $L_+(\Ff)$, such that for all $\nu$ $\in$ $\Ac_{\Ff}$, the process 
$\xi^\nu(\tilde Y)$ is  bounded from below by a martingale.  
The next result gives a characterization of the process $Y$ in terms of backward stochastic differential equation (BSDE)
and of the optimal strategy to problem \reff{vglobal}.

\begin{Thm}  \label{thmbefore}
When  $p$ $>$ $0$ (resp. $p$ $<$ $0$),  the process $Y$ in \reff{defYdynpro} is the smallest (resp. largest) solution in $L_+^b(\Ff)$  to the BSDE:
\beq \label{BSDEY}
Y_t &=& G_T + \int_t^T  f(\theta,Y_\theta,\phi_\theta)  d\theta -  \int_t^T \phi_\theta dW_\theta, \;\;\; 0 \leq t\leq T,
\enq
for some  $\phi$ $\in$  $L_{loc}^2(W)$,  and where
\beq \label{driverf}
f(t,Y_t,\phi_t) &=& p \; \esssup_{\nu \in \Ac_{\Ff}} \Big[ (\mu_t^{\Ff} Y_t + \sigma_t^{\Ff} \phi_t) \nu_t
- \frac{1-p}{2} Y_t |\nu_t\sigma_t^{\Ff}|^2   + K_t^p \frac{(1-\nu_t\gamma_t)^p}{p}  \Big].
\enq
The optimal strategy $(\hat\pi_t^{\Ff})_{t\in [0,T]}$ to problem \reff{vglobal} attains the supremum in \reff{driverf}. Moreover, under the integrability condition:
$\int_0^T \Big|\frac{K_t}{\sigma_t^{\Ff}}\Big|^{\frac{2p}{2-p}}dt$ $<$ $\infty$ a.s.,  the supremum in \reff{driverf} can be taken pointwise, i.e.
\beqs
f(t,Y_t,\phi_t) &=& p \; \esssup_{ \pi < 1/\gamma_t } \Big[ (\mu_t^{\Ff} Y_t + \sigma_t^{\Ff} \phi_t) \pi
- \frac{1-p}{2} Y_t |\pi \sigma_t^{\Ff}|^2   + K_t^p \frac{(1-\pi \gamma_t)^p}{p}  \Big],
\enqs
while  the optimal strategy is given by:
\beqs
\hat\pi_t^{\Ff} &=& \argmax_{\pi < 1/\gamma_t} \Big[ (\mu_t^{\Ff} Y_t + \sigma_t^{\Ff} \phi_t) \pi
- \frac{1-p}{2} Y_t |\pi \sigma_t^{\Ff}|^2   + K_t^p \frac{(1-\pi \gamma_t)^p}{p}  \Big], \;\;\; 0\leq t\leq T, 
\enqs
and satisfies the estimates:
\beq \label{estimpi}
\pi_t^M - \Big(\frac{ \gamma_t^p K_t^p}{(1-p)Y_t |\sigma_t^{\Ff}|^2}\Big)^{\frac{1}{2-p}} & \leq & \hat\pi_t^{\Ff}  \; \leq \; \pi_t^M, \;\;\;\; 0 \leq t\leq T,
\enq
where
\beqs
\pi_t^M &=& \min \Big( \frac{\mu_t^{\Ff}}{(1-p)|\sigma_t^{\Ff}|^2} + \frac{\phi_t}{(1-p)Y_t\sigma_t^{\Ff}} \; , \; \frac{1}{\gamma_t} \Big).
\enqs
\end{Thm}
{\bf Proof.}  By Lemma \ref{lempre}, the process $Y$ lies in $L_+^b(\Ff)$.  
From \reff{dynprosurmar},  we know that for any $\nu$ $\in$ $\Ac_{\Ff}$, the process $\xi^\nu(Y)$ is a $(\Pp,\Ff)$-supermartingale.  
In particular, by taking $\nu$ $=$ $0$, we see that the process $\{Y_t+\int_0^t K_\theta^p d\theta,0\leq t\leq T\}$ is a
$(\Pp,\Ff)$-supermartingale.  From the Doob-Meyer decomposition, and the (local) martingale representation theorem for brownian motion filtration,
we get the existence of  $\phi$ $\in$ $L_{loc}^2(W)$, and a finite variation $\Ff$-adapted process $A$ such that:
\beq \label{dynY}
dY_t  &=& \phi_t dW_t - dA_t, \;\;\; 0 \leq t\leq T.
\enq
From \reff{dynXF} and It\^o's formula, we deduce that the finite variation process in the decomposition of the $(\Pp,\Ff)$-supermartingale $\xi^\nu(Y)$,
$\nu$ $\in$ $\Ac_{\Ff}$, is given by $-A^\nu$ with
\beqs
dA_t^\nu &=& (X_t^{\nu,\Ff})^p \Big\{  \frac{1}{p} dA_t - \Big[  (\mu_t^{\Ff}Y_t + \sigma_t^{\Ff}\phi_t) \nu_t  - \frac{1-p}{2} Y_t |\nu_t\sigma_t^{\Ff}|^2
+ K_t^p \frac{(1-\nu_t\gamma_t)^p}{p}  \Big]  dt \Big\}.
\enqs
Now, by the supermartingale property of $\xi^\nu(Y)$,  $\nu$ $\in$ $\Ac_{\Ff}$, which means that $A^\nu$ is nondecreasing,
and the martingale property of $\xi^{\hat\pi^{\Ff}}(Y)$, i.e. $A^{\hat\pi^{\Ff}}$ $=$ $0$, this implies:
\beqs
 dA_t &=&  p \Big[  (\mu_t^{\Ff}Y_t + \sigma_t^{\Ff}\phi_t) \hat\pi_t^{\Ff}  - \frac{1-p}{2} Y_t |\hat\pi_t^{\Ff}\sigma_t^{\Ff}|^2
+ K_t^p \frac{(1-\hat\pi_t^{\Ff}\gamma_t)^p}{p}  \Big]  dt \\
&=& p \; \esssup_{\nu\in\Ac_{\Ff}}  \Big[  (\mu_t^{\Ff}Y_t + \sigma_t^{\Ff}\phi_t) \nu_t  - \frac{1-p}{2} Y_t |\nu_t\sigma_t^{\Ff}|^2
+ K_t^p \frac{(1-\nu_t\gamma_t)^p}{p}  \Big]  dt.
\enqs
Observing from \reff{defYdynpro} that $Y_T$ $=$ $G_T$,  this proves together with \reff{dynY} that $(Y,\phi)$ solves the BSDE \reff{BSDEY}. In particular,
the process $Y$ is continuous.

Consider now another solution $(\tilde Y,\tilde\phi)$ $\in$ $L_+^b(\Ff)\times L_{loc}^2(W)$ to the BSDE \reff{BSDEY}, and define the family of
nonnegative $\Ff$-adapted processes $\tilde\xi^\nu(\tilde Y)$, $\nu$ $\in$ $\Ac_{\Ff}$, by:
\beq \label{deftildexi}
\xi_t^\nu(\tilde Y) &=&  U(X_t^{\nu,\Ff}) \tilde Y_t + \int_0^t U(X_\theta^{\nu,\Ff}) (1-\nu_\theta\gamma_\theta)^p K_\theta^p d\theta, \;\;\; 0 \leq t\leq T.
\enq
 By It\^o's formula, we see by the same calculations as above that: $d\xi_t^\nu(\tilde Y)$ $=$ $d\tilde M_t^\nu$ $-$ $d\tilde A_t^\nu$, where
 $\tilde A^\nu$ is a nondecreasing $\Ff$-adapted process, and $\tilde M^\nu$ is a local $(\Pp,\Ff)$-martingale as  a stochastic integral with respect to the brownian motion $W$.  By Fatou's lemma under the condition $\tilde Y$ $\in$ $L_+^b(\Ff)$, this implies that the  process $\xi^\nu(\tilde Y)$ is a 
 $(\Pp,\Ff)$-supermartingale, for any $\nu$ $\in$ $\Ac_{\Ff}$.  Recalling that $\tilde Y_T$ $=$ $G_T$, we deduce that for all $\nu$ $\in$ $\Ac_{\Ff}$
 \beqs
 \Ee \Big[ U(X_T^{\nu,\Ff}) G_T + \int_t^T  U(X_\theta^{\nu,\Ff}) (1-\nu_\theta\gamma_\theta)^p K_\theta^p d\theta \Big| \Fc_t \Big] & \leq&
 U(X_t^{\nu,\Ff}) \tilde Y_t,  \;\;\; 0 \leq t\leq T. 
 \enqs
If $p$ $>$ $0$ (resp. $p$ $<$ $0$), then by dividing the above inequalities by  $U(X_t^{\nu,\Ff})$, which is positive (resp. negative), we deduce 
by definition of $Y$ (see \reff{Yppos} and \reff{Ypneginter}), and arbitrariness of $\nu$ $\in$ $\Ac_{\Ff}$,  that   $Y_t$ $\leq$ (resp. $\geq$)  $\tilde Y_t$, $0\leq t\leq T$. 
This shows that $Y$ is the smallest (resp. largest) solution to the BSDE \reff{BSDEY}. 

Next, we make the additional integrability condition:
\beq \label{integK}
\int_0^T \Big|\frac{K_t}{\sigma_t^{\Ff}}\Big|^{\frac{2p}{2-p}}dt &<& \infty, \;\;\; a.s.
\enq
Let us consider the function $F$ defined on $\{(\omega,t,\pi) \in \Omega\times [0,T]\times \R: \pi < 1/\gamma_t(\omega)\}$ by:
\beqs
F(t,\pi) &=&  (\mu_t^{\Ff} Y_t + \sigma_t^{\Ff} \phi_t) \pi
- \frac{1-p}{2} Y_t |\pi \sigma_t^{\Ff}|^2   + K_t^p \frac{(1-\pi \gamma_t)^p}{p}.
\enqs
(As usual, we omit the dependence of $F$ in $\omega$).  By definition, we clearly have almost surely
\beq \label{fpoint}
\frac{1}{p} f(t,Y_t,\phi_t) & \leq & \esssup_{\pi < 1/\gamma_t} F(t,\pi), \;\;\; 0 \leq t\leq T. 
\enq
Let us prove the converse inequality. Observe that, almost surely,  for all $t$ $\in$ $[0,T]$, the function $\pi$ $\rightarrow$
$F(t,\pi)$ is  strictly concave (recall that the process $Y$ is strictly positive), $C^2$ on $(-\infty,1/\gamma_t)$, with:
\beqs
\Dpi{F}(t,\pi) &=&  (\mu_t^{\Ff} Y_t + \sigma_t^{\Ff} \phi_t) - (1-p) Y_t |\sigma_t^{\Ff}|^2 \pi - \gamma_t K_t^p (1-\pi\gamma_t)^{p-1},
\enqs
and satisfies:
\beqs
\lim_{\pi\rightarrow -\infty} F(t,\pi) \; = \; - \infty,  & & \lim_{\pi\rightarrow -\infty} \Dpi{F}(t,\pi) \; = \; \infty, \;\;
\lim_{\pi\rightarrow 1/\gamma_t} \Dpi{F}(t,\pi) \; = \; - \infty. 
\enqs
We deduce that almost surely,   for all $t$ $\in$ $[0,T]$, the function
$\pi$ $\rightarrow$ $F(t,\pi)$ attains  its maximum at some point $\hat\pi_t^{\Ff}$, which satisfies:
\beqs
\Dpi{F}(t,\hat\pi_t^{\Ff}) &=& 0. 
\enqs
By a measurable selection argument,  this defines an $\Ff$-adapted process $\hat\pi^{\Ff}$ $=$ $(\hat\pi_t^{\Ff})_{t\in [0,T]}$.
In order to prove the equality in \reff{fpoint}, it suffices to show that such $\hat\pi^{\Ff}$ lies in $\Ac_{\Ff}$, and this will be checked under the condition
\reff{integK}.  For this, consider the  $\Ff$-adapted processes  $\tilde\pi^M$ and $\pi^M$  defined by:
\beqs
\tilde\pi_t^M &=&  \frac{\mu_t^{\Ff}}{(1-p)|\sigma_t^{\Ff}|^2} + \frac{\phi_t}{(1-p)Y_t\sigma_t^{\Ff}}, \;\;\;
\pi_t^M \; = \; \min\big(\tilde \pi_t^M, \frac{1}{\gamma_t}\big),  \;\;\; 0 \leq t\leq T.
\enqs
When $\tilde\pi_t^M(\omega)$ $<$ $1/\gamma_t(\omega)$, we have:
\beqs
\Dpi{F}(t,\tilde\pi_t^M) &=& - \gamma_t K_t^p (1-\tilde\pi_t^M\gamma_t)^{p-1} \; \leq \; 0,
\enqs
and so by strict concavity of $F(t,\pi)$ in $\pi$: $\hat\pi_t^{\Ff}$ $\leq$ $\tilde\pi_t^M$. When $\tilde\pi_t^M(\omega)$ $\geq$ $1/\gamma_t(\omega)$,
and since
$\hat\pi_t^{\Ff}$ $<$ $1/\gamma_t$, we get:   $\hat\pi_t^{\Ff}$ $\leq$ $\tilde\pi_t^M$.  Consequently, we have the upperbound:
$\hat\pi_t^{\Ff}$ $\leq$ $\pi_t^M$, for all $t$ $\in$ $[0,T]$.
Notice that by \reff{integmu}, continuity of the path of $Y$, and since $\phi$ $\in$  $L_{loc}^2(W)$, we have:
$\int_0^T |\tilde\pi_t^M\sigma_t^{\Ff}|^2 dt$ $<$ $\infty$ a.s. Moreover, since $\gamma_t$ $\geq$ $0$, we have $|\pi^M|$ $\leq$ $|\tilde\pi^M|$, and thus
$\pi^M$ lies in $\Ac_{\Ff}$.  Next, consider the $\Ff$-adapted process $\bar\pi$ defined by:
\beqs
\bar\pi_t &=& \pi_t^M - \rho_t, \;\;\; 0 \leq t \leq T,
\enqs
for some $\Ff$-adapted nonnegative process $\rho$ $=$ $(\rho_t)_{t\in [0,T]}$ to be determined. When $\tilde\pi_t^M(\omega)$ $<$
$1/\gamma_t(\omega)$, we have
\beq
\Dpi{F}(t,\bar\pi_t) &=&  (1-p)Y_t|\sigma_t^{\Ff}|^2\rho_t - \gamma_t K_t^p (1-\tilde\pi_t^M\gamma_t + \rho_t\gamma_t)^{p-1}  \nonumber \\
& \geq &  (1-p)Y_t|\sigma_t^{\Ff}|^2\rho_t - \gamma_t K_t^p (\rho_t\gamma_t)^{p-1}.  \label{Dfpiineg}
\enq
When $\tilde\pi_t^M(\omega)$ $\geq$ $1/\gamma_t(\omega)$, the inequality \reff{Dfpiineg} also holds true. Hence, by choosing $\rho$ such that the r.h.s. of \reff{Dfpiineg} vanishes, i.e.
\beqs
\rho_t &=& \Big(\frac{ \gamma_t^p K_t^p}{(1-p)Y_t |\sigma_t^{\Ff}|^2}\Big)^{\frac{1}{2-p}}, \;\;\; 0 \leq t\leq T,
\enqs
we obtain almost surely: $\Dpi{F}(t,\bar\pi_t)$ $\geq$ $0$, $0\leq t\leq T$, and so by strict concavity of $F$ in $\pi$: $\bar\pi_t$ $\leq$ $\hat\pi_t^{\Ff}$.
Finally, under \reff{integK}, and recalling that $Y$ is continuous, $\gamma$ $<$ $1$,  we easily see that $\rho$ satisfies the integrability condition:
$\int_0^T |\rho_t\sigma_t^{\Ff}|^2 dt$ $<$ $\infty$ a.s., and so  $\bar\pi$ lies in $\Ac_{\Ff}$. Therefore, we have proved that $\hat\pi^{\Ff}$ lies in $\Ac_{\Ff}$, and satisfies the estimates \reff{estimpi}.
\ep

\begin{Rem} \label{rempneg}
{\rm  The driver $f(t,Y_t,\phi_t)$ of the BSDE \reff{BSDEY} is in general not Lipschitz in the arguments in $(Y_t,\phi_t)$,  and we are not able to prove  by standard arguments that there exists a unique solution to this BSDE. 
}
\end{Rem}

\begin{Rem}
{\rm  We make some comments and interpretation on the form of the optimal before-default strategy.  Let us consider a default-free stock
market model  with drift and volatility coefficients $\mu^{\Ff}$ and $\sigma^{\Ff}$,
and an investor with CRRA utility function $U(x)$ $=$ $x^p/p$,   looking for the optimal investment problem:
\beqs
V_0^{M} &=& \sup_{\pi\in\Ac_{\Ff}} \Ee[ U(X_T^{\Ff})],
\enqs
where $X^{\Ff}$  is the wealth process in \reff{dynXF}.
In this context, $\Ac_{\Ff}$, defined in \reff{defAF} is interpreted as the set of trading strategies that are constrained to be upper-bounded (in proportion)
by  $1/\gamma_t$.  In other words, $V_0^{M}$ is the Merton optimal investment problem under constrained strategies.  By considering, similarly as in
\reff{defYdynpro}, the process
\beqs
Y_t^{M} &=&  p \esssup_{\pi^{\Ff} \in \Ac_{\Ff}(t,\nu)} \Ee  \Big[ U\Big(\frac{X_T^{\Ff}}{X_t^{\nu,\Ff}}\Big) \Big| \Fc_t \Big], \;\;\;  0 \leq t\leq T,
\enqs
and arguing similarly as in Theorem \ref{thmbefore}, one can prove that $Y^{M}$ is the smallest solution  to the BSDE:
\beqs
Y_t^M &=& 1 + \int_t^T  f^M(\theta,Y_\theta^M,\phi_\theta^M)  d\theta -  \int_t^T \phi_\theta^M dW_\theta, \;\;\; 0 \leq t\leq T,
\enqs
for some  $\phi^M$ $\in$  $L_{loc}^2(W)$,  where
\beqs
f^M(t,Y_t^M,\phi_t^M) &=& p \; \esssup_{ \pi < 1/\gamma_t } \Big[ (\mu_t^{\Ff} Y_t^M + \sigma_t^{\Ff} \phi_t^M) \pi
- \frac{1-p}{2} Y_t^M |\pi \sigma_t^{\Ff}|^2  \Big],
\enqs
while the optimal strategy for $V_0^M$ is given by:
\beqs
\hat\pi_t^{M} &=& \min \Big( \frac{\mu_t^{\Ff}}{(1-p)|\sigma_t^{\Ff}|^2} + \frac{\phi_t^M}{(1-p)Y_t^M\sigma_t^{\Ff}} \; , \; \frac{1}{\gamma_t} \Big).
\enqs
Notice that when the coefficients $\mu^{\Ff}$, $\sigma^{\Ff}$ and $\gamma$ are deterministic, then $Y^M$ is also deterministic, i.e. $\phi^M$ $=$ $0$, and is the positive solution to the ordinary differential equation:
\beqs
Y_t^M &=& 1 + \int_t^T  f^M(\theta,Y_\theta^M)  d\theta, \;\;\; 0 \leq t\leq T,
\enqs
with $f^M(t,y)$ $=$ $py \sup_{\pi< 1/\gamma_t}[\mu_t^{\Ff}\pi - \frac{1-p}{2}|\pi\sigma_t^{\Ff}|^2]$ $=:$ $pyc(t)$, and so: $Y_t^M$ $=$
$\exp(p \int_t^T c(\theta)d\theta)$.
Moreover, the optimal strategy is $\hat\pi_t^M$ $=$ $ \min\big( \frac{\mu_t^{\Ff}}{(1-p)|\sigma_t^{\Ff}|^2}\; , \; \frac{1}{\gamma_t} \big)$.
In particular, when there is no constraint on trading strategies, i.e.
$\gamma$ $=$ $0$, we recover the usual expression of the optimal Merton trading strategy: $\hat\pi_t^M$ $=$ $\frac{\mu_t^{\Ff}}{(1-p)|\sigma_t^{\Ff}|^2}$.

Here, in our default stock market model, the optimal  before-default  strategy $\hat\pi^{\Ff}$ satisfies the estimates \reff{estimpi}, which have the following interpretation. The process  $\pi^M$ has a similar form as the optimal Merton trading strategy $\hat\pi^M$ described above, but includes further
through the process $Y$ and $K$, the eventuality of a default of the stock price, inducing a drop of size $\gamma$, and then
a switch of the coefficients of the stock price from $(\mu^{\Ff},\sigma^{\Ff})$ to $(\mu^d,\sigma^d)$.  The optimal trading strategy $\hat\pi^{\Ff}$ is
upper-bounded by  $\pi^M$, and when the jump size $\gamma$ goes to zero, it converges to $\hat\pi^M$, as expected since  in this case the
model behaves as a no-default market.

}
\end{Rem}

\subsection{Example and numerical illustrations}

We consider a special case where $\mu^{\Ff}$,
$\sigma^{\Ff}$, $\gamma$ are constants, $\mu^d(\theta)$
$\sigma^d(\theta)$ are only deterministic functions of
$\theta$, and the default time $\tau$ is independent of
$\Ff$, so that  $\alpha_t(\theta)$ is simply  a known
deterministic function $\alpha(\theta)$ of  $\theta$
$\in$ $\R_+$, and the survival probability   $G(t)$ $=$
$\Pp[\tau > t |\Fc_t]$ $=$ $\Pp[\tau> t]$ $=$
$\int_t^\infty \alpha(\theta) d\theta$ is a
deterministic function. We also choose a CRRA utility
function $U(x)$ $=$ $\frac{x^p}{p}$, $p$ $<$ $1$,
$p\neq 0$, $x$ $>$ $0$.  Notice that $V_\theta^d(x)$
$=$ $v^d(\theta,x)$ $=$ $U(x)k(\theta)^p$ with 
\beqs
k(\theta) &=&  \Big( \Ee\Big[\alpha_T(\theta)
\Big(\frac{Z_T(\theta)}{\alpha_T(\theta)} \Big)^{-q}
\Big] \Big)^{1\over q} \; = \; \alpha(\theta)^{\frac
1p} \exp\Big( \frac{1}{2}
\Big|\frac{\mu^d(\theta)}{\sigma^d(\theta)}\Big|^2
\frac{1}{1-p} (T-\theta) \Big) 
\enqs 
Moreover,  the optimal wealth process after-default does not depend on
the default time density, and the optimal strategy
after-default is  given, similarly as in the Merton
case, by: \beqs \hat\pi^d(\theta) &=&
\frac{\mu^d(\theta)}{(1-p)|\sigma^d(\theta)|^2}. \enqs
On the other hand, from the above results and
discussion, we know that in this Markovian case, the value function of the global before-default optimization
problem is in the form $V_0$ $=$ $v(0,X_0)$ with: 
\beqs
v(t,x) &=& U(x) Y(t), 
\enqs 
where $Y$ is a deterministic function of time, solution to the first-order ordinary differential equation (ODE):
\beq \label{equ:ode numerical example} 
Y(t) &=& G(T) + \int_t^T f(\theta,Y(\theta)) d\theta, \;\;\; t \in [0,T],
\enq
with 
\beq \label{fpi}
f(t,y) &=& p \sup_{\pi< 1/\gamma} \Big[ \Big(\mu^{\Ff} \pi  - \frac{1-p}{2} |\pi\sigma^{\Ff}|^2\Big) y  + k(t)^p\frac{(1-\pi\gamma)^p}{p} \Big]
\enq
There is no explicit solutions to this ODE, and we shall give some numerical illustrations.

The following numerical results are based on the  model parameters described  below. We suppose
that the survival probability follows the exponential 
distribution with constant default intensity, i.e.
$G(t)=e^{-\lambda t}$ where $\lambda >0$, and thus  the density function is $\alpha(\theta)=\lambda e^{-\lambda \theta}$. 
The functions $\mu^d(\theta)$ and $\sigma^d(\theta)$ are supposed  to be in the form
\beqs
\mu^d(\theta) \; = \; \mu^{\mathbb F}\frac{\theta}{T}, & &
\sigma^d(\theta) \; = \; \sigma^{\mathbb F}(2-\frac{\theta}{T}), \;\;\; \theta \in [0,T],
\enqs
which have  the following economic interpretation. The ratio between the after and before-default rate of return is smaller than one, and increases linearly 
with the default time:  the after-default rate of return drops to zero, when the default time occurs near the initial date, and 
converges  to the  before-default rate of return, when  the default time occurs near the finite investment horizon.  
We have a similar interpreation for the volatility but with symmetric  relation: 
the ratio between the  after and before-default volatility is larger than one,  decreases linearly with the default time, converging to the double (resp. initial) value of the before-default  volatility, when the default time goes to the initial (resp. terminal horizon) time. 
 
To solve numerically the  ODE  \eqref{equ:ode numerical example}, we  apply the Howard algorithm,  which consists in 
iterating in \reff{fpi} the control value $\pi$   at each step of  the ODE  resolution. We initialize the algorithm by choosing the constrained 
Merton strategy
\beqs
\hat \pi^{M} &=& \min\Big(\frac{\mu^{\mathbb F}}{(1-p)|\sigma^{\mathbb F}|^2},\,\frac{1}{\gamma}\Big). 
\enqs 
In the following Table \ref{Tab: strategy gamma}, we
show the impact of the loss given default $\gamma$ on
the optimal strategy $\hat\pi^{\mathbb F}_t$. The
numerical tests show that except in some extreme cases
where both the default probability and the loss given
default are large, the optimal strategy  is quite
invariant with respect to time $t$ in most cases we
consider. So we give below the optimal strategy as its
expected  value on time. We perform numerical results for various  degrees of risk aversion $1-p$:  smaller, close to  and larger than one, 
and with $\mu^{\mathbb F}=0.03$, $\sigma^{\mathbb F}=0.1$, $T$ $=$ $1$ and $\lambda=0.01$.

\begin{table}[h]
\caption{Optimal strategy vs constrained Merton. 
}
\label{Tab: strategy gamma}
\begin{center}
\begin{tabular}{|c||c|c|c|c|c|c|}\hline
&\multicolumn{2}{c|}{$p=0.2$}&\multicolumn{2}{c|}{$p\rightarrow 0$}
&\multicolumn{2}{c|}{$p=-0.2$}\\
\hline\parbox{2mm}{$\gamma$}
&\parbox{2mm}{$\hat\pi^{\mathbb
F}$}&\parbox{6mm}{$\hat\pi^{\mathrm{M}}$}&\parbox{2mm}{$\hat\pi^{\mathbb
F}$}&\parbox{6mm}{$\hat\pi^{\mathrm{M}}$}&\parbox{2mm}{$\hat\pi^{\mathbb
F}$}&\parbox{6mm}{$\hat\pi^{\mathrm{M}}$}\\\hline
$0.01$&$3.73$&$3.74$&$2.99$&$3.00$&$2.49$&$2.50$\\
$0.1$&$3.57$&$3.74$&$2.86$&$3.00$&$2.38$&$2.50$\\
$0.5$&$1.58$&$2.00$&$1.38$&$2.00$&$1.22$&$2.00$\\
$0.8$&$0.91$&$1.25$&$0.80$&$1.25$&$0.70$&$1.25$\\
\hline
\end{tabular}
\end{center}
\end{table}

First, observe that, when we take into account the 
counterparty risk, the  proportion  invested in the stock
is always smaller than the Merton strategy without
conterparty risk.  Secondly, the  strategy is decreasing with respect to $\gamma$, which
means that one should reduce the stock investment when
the loss given default increases. 

\begin{table}[h]\caption{Optimal strategy with various $\lambda$
and $\gamma$.}\label{Tab:strategy lambda}
\begin{center}
\begin{tabular}{|c|c||c|c|c|c|c|c|}\hhline{|--|------}
\multicolumn{2}{|c||}{}&\multicolumn{2}{c|}{$p=0.2$}&\multicolumn{2}{c|}{$p\rightarrow
0$} &\multicolumn{2}{c|}{$p=-0.2$}\\\hhline{|~~|------}
\multicolumn{2}{|c||}{}&$\gamma=0.1$&$\gamma=0.5$&$\gamma=0.1$&$\gamma=0.5$&$\gamma=0.1$&$\gamma=0.5$
\\\hline
\multicolumn{2}{|c||}{\parbox{6mm}{$\hat\pi^{\mathrm{M}}$}}
&$3.75$&$2.00$&$3.00$&$2.00$&$2.50$&$2.00$\\\hhline{|==#=|=|=|=|=|=|}
$\lambda=0.01$&PD$=0.01$&$3.57$&$1.58$&$2.86$&$1.38$&$2.38$&$1.22$\\\hhline{}
$\lambda=0.05$&PD$=0.05$&$2.93$&$0.26$&$2.35$&$0.21$&$1.96$&$0.18$\\
$\lambda=0.1$&PD$=0.10$&$2.22$&$-0.90$&$1.78$&$-0.70$&$1.49$&$-0.58$\\
$\lambda=0.3$&PD$=0.26$&$0.00$&$-3.96$&$0.00$&$-3.00$&$0.00$&$-2.40$\\
\hline
\end{tabular}
\end{center}
\end{table}

Next, we examine  the role played by the default intensity $\lambda$. In
Table \ref{Tab:strategy lambda}, the column PD represents the
default probability of the counterparty up to $T$ with
the given intensity $\lambda$, i.e. PD  $=$ $\proba(\tau\leq T)=1-e^{-\lambda T}$. We see that the stock
investment decreases rapidly when the default
probability increases. Moreover, when both default
probability of the counterparty and the loss given
default of the stock are large, one should take short
position on the stock in the portfolio investment
strategy before the default of the counterparty. Then at the default time $\tau=\theta$, 
the optimal strategy is switched to $\hat\pi^d(\theta)$, which is always positive.

\begin{figure}[h]\caption{Value function -- optimal vs Merton: $p=0.1$, 
$\lambda=0.01$.}\label{fig:value function gamma}
\begin{center}
\epsfig{file=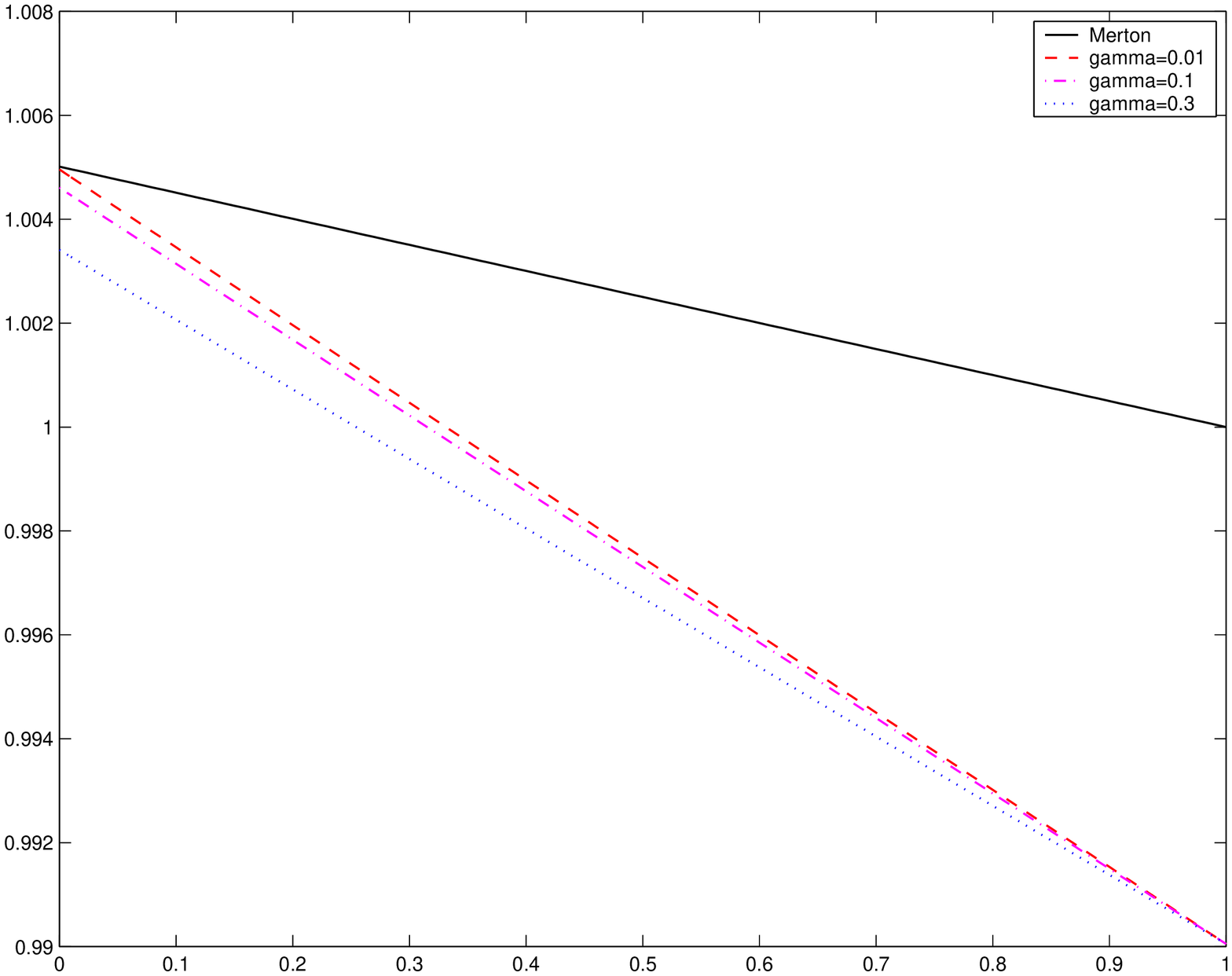,height=9cm,angle=90} 
\end{center}
\caption{Value function -- optimal vs Merton: $p=0.1$,
$\gamma=0.1$.}\label{fig:value function lambda}
\begin{center}
\epsfig{file=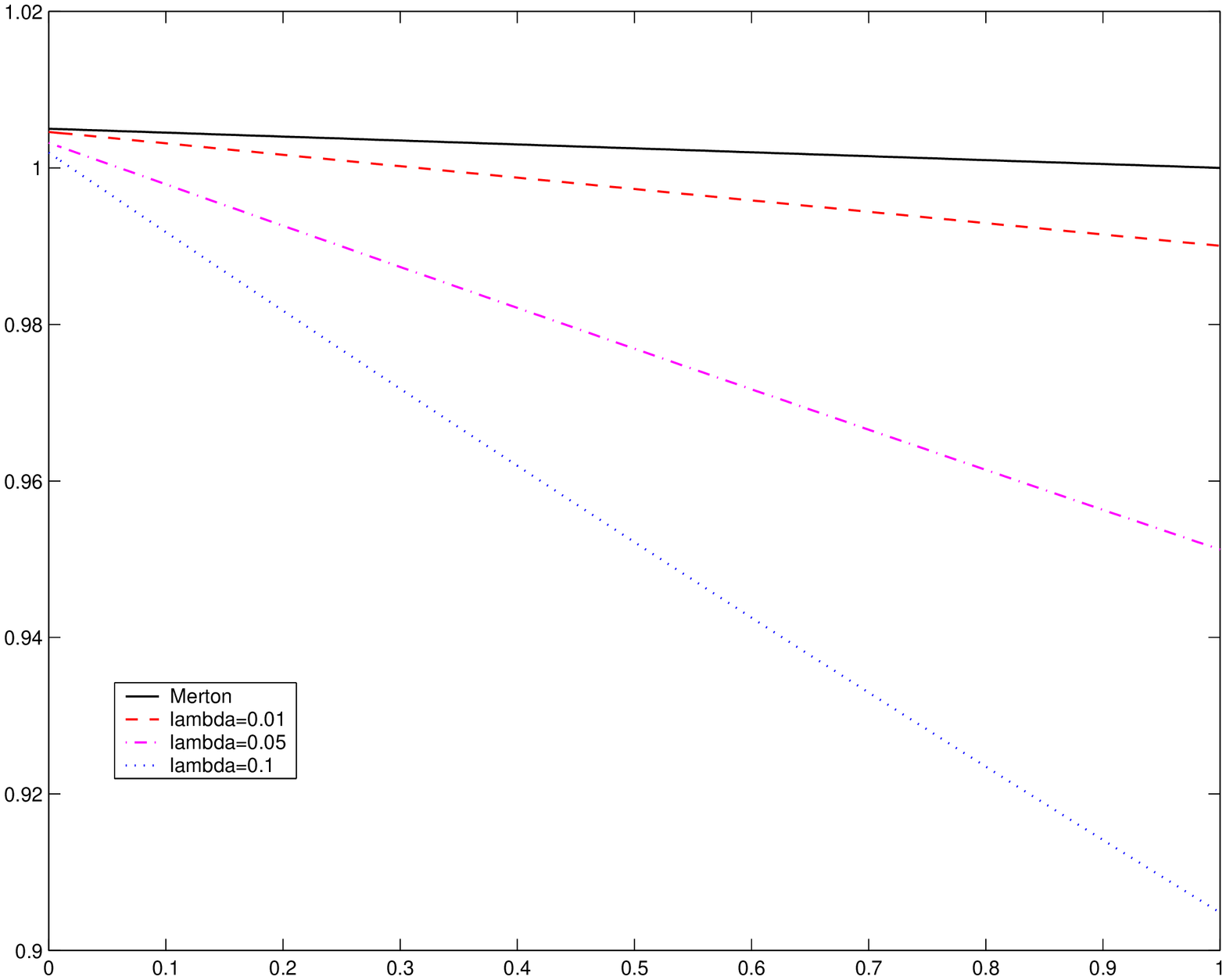,height=9cm,angle=90}
\end{center}
\end{figure}

We also compare the value function obtained in our
example to that in the classical Merton model, that is,
the solution $Y(t)$ to the ODE \eqref{equ:ode numerical
example} and the function $Y^{M}(t)$ deduced with
$k(t)=0$ and $G(T)=1$. In Figure \ref{fig:value
function gamma}, the curves represent different values of $\gamma$ such that 
$\hat \pi^{M}$ $=$ $\frac{\mu^{\mathbb F}}{(1-p)|\sigma^{\mathbb F}|^2}\leq \frac{1}{\gamma}$. 
The value function $Y(t)$ obtained with counterparty risk is always
below the Merton one. It is decreasing on time and  also decreasing w.r.t. the 
proportional loss $\gamma$. In addition, for a given default intensity $\lambda$,
all curves converge at $T$ to $G(T)=e^{-\lambda T}$. In
Figure \ref{fig:value function lambda}, the curves
represent different values of $\lambda$. The value
function $Y(t)$ is also decreasing w.r.t. the default intensity
$\lambda$. However, the final value of each curve corresponds to different values of $G_T$.

\clearpage

\clearpage


\section{Conclusion}

This paper studies an optimal investment problem under the presence of counterparty risk for the trading stock. By adopting a conditional density approach for the default time, we derive  a suitable decomposition of the initial utility maximization problem into an after-default  and a global default one, the solution to the latter depending on the former.  This makes  the resolution of the optimization problem more explicit, and provides a fine understanding 
of the optimal trading  strategy emphasizing the impact of default time and loss given default.  
The density approach  can be used for studying other  optimal  portfolio problems,   
like  the pricing by indifference-utility,  with counterparty risk.  A further topic is the  optimal investment problem  with two assets (names) exposed 
both to  bilateral counterparty risk, and the  conditional density approach  should be relevant  for such study  planned for future research.


\end{document}